\documentclass[12pt,a4paper]{amsart}
\usepackage{amsmath, amssymb,mathrsfs}

\numberwithin{equation}{section}
\newcommand{\wedg}{\mathbin{\scriptstyle{\wedge}}}
\newcommand{\lctr}{\mathbin{\lrcorner}}
\newcommand{\contr}{\lctr\,}
\newcommand{\Wedge}{\Lambda}

\newcommand{\supp}{\operatorname{supp}}
\newcommand {\rot}{\operatorname{\mathbf{curl}}}

\newcommand {\grad}{\operatorname{\mathbf{grad}}}
\newcommand {\Div}{\operatorname{div}}

\newcommand{\Cinf}{{\mathscr C}^\infty}
\let\ov\overline
\newcommand{\C}{\mathbb{C}}
\newcommand{\R}{\mathbb{R}}
\newcommand{\N}{\mathbb{N}} 
 
\newcommand{\Z}{\mathbb{Z}}
\newcommand{\Adot}{\dot{\mathscr A}}
\newcommand{\Bdot}{\dot{B}}
\newcommand{\Fdot}{\dot{F}}
\newcommand{\Hdot}{\dot{H}}

\newcommand{\J}{\mathscr{J}}
\newcommand{\tent}{\mathcal T}
\newcommand{\ap}{\sigma}
\newcommand{\app}{\beta}
\newcommand{\ind}{\mathbf 1}
\newcommand{\norm}[1]{\left\Vert #1 \right\Vert}
\newcommand{\ds}{\displaystyle}

\newtheorem{theorem}{Theorem}[section]
\newtheorem{lemma}[theorem]{Lemma}  
\newtheorem{proposition}[theorem]{Proposition}
\newtheorem{corollary}[theorem]{Corollary}
\newtheorem{remark}[theorem]{Remark}

\newtheorem{definition}[theorem]{Definition}

\newenvironment{proof1}{\parindent=0pt \parskip=0.5ex 
   \abovedisplayskip = 0.8 \abovedisplayskip
   \belowdisplayskip=\abovedisplayskip
   \textbf{Proof: \ }}{\hfill $\square$\smallskip}
\usepackage{color}
\title[Potential maps, Hardy spaces, and tent spaces on domains]{Potential maps, Hardy spaces, and tent spaces on special Lipschitz domains}

\author{Martin Costabel \and Alan McIntosh \and Robert J. Taggart}

\date{20 February, 2012}

\begin{document}

\maketitle

\begin{abstract}

\noindent Suppose that $\Omega$ is the open region in $\R^n$ above a Lipschitz graph and let $d$ denote the exterior derivative on $\R^n$. We construct a convolution operator $T $ which preserves support in $\ov\Omega$, is smoothing of order $1$ on the homogeneous function spaces, and is a potential map in the sense that $dT$ is the identity on spaces of exact forms with support in $\ov\Omega$. Thus if $f$ is exact and supported in $\ov\Omega$, then there is a potential $u$, given by $u=Tf$, of optimal regularity and supported in $\ov\Omega$, such that $du=f$. This has implications for the regularity in homogeneous function spaces of the de Rham complex on $\Omega$ with or without boundary conditions. The operator $T$ is used to obtain an atomic characterisation of Hardy spaces $H^p$ of exact forms with support in $\ov\Omega$ when $n/(n+1)<p\leq1$.   This is done via an atomic decomposition of functions in the  tent spaces $\tent^p(\R^n\times\R^+)$  with support in a tent $T(\Omega)$ as a sum of atoms with support  away from the boundary of $\Omega$.
This new decomposition of tent spaces is useful, even for scalar valued functions.

\medskip



\end{abstract}

\section{Introduction}\label{s:intro}

The study of potential maps on domains $\Omega$ of $\R^n$ has a rich history. Consider, for a moment, the following question. Suppose that a function $f$ belongs to a Sobolev space $H^m_0(\Omega)$ where $m\geq0$ and $\Omega$ is a bounded strongly Lipschitz domain, and suppose that $\int f=0$. Is there a 
vector field $u$ in  $(H^{m+1}_0(\Omega))^n$ which satisfies $\Div u=f$? The answer is yes, as was essentially proved by J.~Ne{\v{c}}as. Indeed, this follows by duality from \cite[Chapter 3, Lemma 7.1]{N}. An alternative proof
was provided by M.~E.~Bogovski\u{\i} \cite{mB79,mB80}, who, in particular, constructed an integral operator $T$ which maps boundedly from the Sobolev space $W^{m,p}_0(\Omega)$ into $(W^{m+1,p}_0(\Omega))^n$ in the case when $m\geq 0$, $1<p<\infty$, and $\Omega$ is starlike with respect to a ball. The potential $u$ that solves the equation $\Div u=f$ is then given by $u=Tf$, provided that $\int f=0$.  Since $T$ gives the potential $u$ which solves the equation and preserves support, we say that $T$ is a {\it potential map} for the domain $\Omega$.  Such a potential map is an important tool in the theory of equations of hydrodynamics. See the monograph \cite{gG94} of G.~P.~Galdi and the papers mentioned below for further references to the extensive literature.

Subsequently M.~Mitrea \cite{mM04} and D.~Mitrea, M.~Mitrea and S.~Monniaux \cite{MMM08} 
adapted Bogovski\u{\i}'s operator to construct potential maps $T_{\ell}$ to solve the equation $du=f$, where $d$ denotes the exterior derivative operator and where $f$ is an exact $\ell$-form with support in $\ov\Omega$ and coefficients in a suitable function space. Using $T_{\ell}$, they have thereby deduced sharp regularity estimates for important operators in the theory of hydrodynamics. As with Bogovski\u{\i}'s operator, the mapping properties of $T_{\ell}$ ensure that there is no loss of regularity, and that support in $\ov\Omega$ is preserved. In the case when $n=3$, the equation $du=f$ is equivalent to solving one of the equations 
\[\grad u=f, \qquad \rot u =f \qquad\mbox{or}\qquad \Div u=f,\]
where $f$ is interpreted either as a scalar or vector field depending on the value of $\ell$.

More recently, M.~Costabel and A.~McIntosh~\cite{CoMcI} showed that the operators $T_{\ell}$ are pseudodifferential operators of order $-1$ and are therefore bounded in all the spaces of $\ell$-forms with coefficients in one of the Besov or Triebel--Lizorkin classes.  
For a domain starlike with respect to a ball, the special support properties of the operators imply regularity for the de Rham complex  with full Dirichlet boundary conditions. Similar results hold for complexes without boundary conditions (using dual Poincar\'e-type operators). For bounded strongly Lipschitz domains, the same regularity results hold, and in addition the cohomology spaces can always be represented by $\Cinf$ functions.

In the present paper, we turn our attention to unbounded special Lipschitz domains, that is, to domains in $\R^n$ that lie above the graph of a Lipschitz function, and consider the spaces $\Hdot^s(\R^n,\Wedge)$ of forms with coefficients in the homogeneous Sobolev space of degree $s$, where $s\in\R$.
 We construct an operator $T$ with the following properties. First, $T$ boundedly lifts forms in $\Hdot^s(\R^n,\Wedge)$ to forms in $\Hdot^{s+1}(\R^n,\Wedge)$. Second, $T$ preserves support in $\ov\Omega$. Third, if $df=0$ then a solution $u$ of the equation $du=f$ is given by $u=Tf$. Hence the equation $du=f$ is solved on $\ov\Omega$ with optimal regularity, because if $f\in\dot H^s(\R^n,\Wedge^\ell)$, $df=0$, and the support of $f$ is in $\ov\Omega$, then $u\in \dot H^{s+1}(\R^n,\Wedge^{\ell-1})$, $du=f$, and the support of $u$ is in $\ov\Omega$.

Our potential map $T$ is a convolution integral operator  involving kernels which  are supported in a cone of $\R^n$ lying above its base at the origin. This support property is ideally suited for boundary problems on a given special Lipschitz domains $\Omega$, since such a cone (with appropriately chosen aperture) is contained within $\Omega$ when translated by points from $\Omega$, and thus $T$ preserves support in $\ov\Omega$. The bounded lifting property of $T$ is valid not only when $T$ acts on $\Hdot^s(\R^n,\Wedge)$, but also when $T$ acts on the space of forms whose coefficients belong to any homogeneous Besov or Triebel--Lizorkin function space. These properties of $T$ have immediate ramifications for the regularity of the de Rham complex in such spaces on $\ov\Omega$.

Our formula for $T$ is motivated by a reproducing formula of D.~Chang, S.~Krantz and E.~Stein \cite{CKS93}, which was also constructed to preserve support in a special Lipschitz domain.  

In the second half of this paper, we show how our methods provide an alternative proof of this reproducing formula, and then use it to show that Hardy spaces $H^p_d(\R^n,\Wedge)$ of exact forms can be characterised by exact atoms  whenever $n/(n+1)<p\leq1$. Using the same method, we also show that the Hardy space $H^p_{z,d}(\ov\Omega,\Wedge)$ of exact forms with support in the closure  $\ov\Omega$ of a special Lipschitz domain can be characterised using exact atoms supported in $\ov\Omega$. These results generalise the classical theorems of R.~Coifman \cite{rC74} and R.~Latter \cite{rL78} for $\R^n$, and the result of Chang, Krantz and Stein \cite{CKS93} for special Lipschitz domains, and have already been proved by Z.~Lou and A.~McIntosh \cite{LM05, LM04} for exact forms in the case when $p=1$. In the three latter papers, the authors follow the method of R.~Coifman, Y.~Meyer and E.~Stein \cite{CMS85} to show that each element $u$ of the Hardy space can be written as a sum of atoms by first mapping $u$ into a `tent space' via an operator $Q$, then decomposing the image $Qu$ as a sum of tent space atoms, before finally mapping the decomposition back into the Hardy space. In the case that $u$ is supported in $\ov\Omega$, \cite{CKS93} and \cite{LM04} then use reflection maps to express $u$ as a sum of atoms whose supports lie in $\ov\Omega$. Our method of proof differs from both of these papers, in that we use the full power of the reproducing formula to obtain a tent space atomic decomposition of $Qu$ whose atoms have good support properties with respect to $\ov\Omega$, thus removing the need for reflection maps.

 This result on tent spaces is new and interesting, even for scalar valued functions. For this reason, we state here a simplified version of Theorem \ref{th:tent space decomposition}. For details, see Section \ref{s:atomic decomp}.   What is new is the requirement that $5B_k\subset\Omega$.

\begin{theorem} \label{1.1}
Suppose $U$ belongs to the tent space $\tent^1(\R^n\times\R^+)$, with $\supp(U)$ in a tent $T(\Omega) $ over an open set $\Omega\subset\R^n$. Then $U=\sum_k\lambda_kA_k$ where  $\sum|\lambda_k| \lesssim \|U\|_{\tent^1}$ and $\supp (A_k)\subset B_k\times (0,6 r_k)$ where $B_k$ is a ball with radius $r_k$ such that $5 B_k\subset\Omega$, and $\text{vol}(B_k)\int_0^\infty\|A_k\|_2^2\,\frac{dt}t= 1$.
\end{theorem} 

The paper is organised as follows. In Section~\ref{s:notation} we introduce notation and define the various spaces that we use. At the beginning of Section~\ref{s:special lipschitz}, the potential map $T$ is defined, its properties are stated (see Theorem~\ref{th:T on special lipschitz}) and its utility for solving potential equations with boundary conditions on special Lipschitz domains is illustrated. The rest of this section is devoted to the proof of Theorem~\ref{th:T on special lipschitz}. Section~\ref{s:restriction spaces} contains applications for potential maps on spaces without boundary conditions, by considering operators induced from $T$ via quotient maps.  In Section 5, we show how our expressions  relate to the reproducing formula of Chang, Krantz and Stein, then in the final section, we  use it to  prove our atomic characterisation of Hardy spaces of exact forms on $\R^n$ and $\ov\Omega$.

\section*{Acknowledgments}

The authors appreciate the support of the Centre for Mathematics and its Applications at the Australian National University, Canberra, where this project was commenced during a visit of the first author. Support from the Australian Government through the Australian Research Council is also gratefully acknowledged.

The authors would like to thank Andrew Morris for helpful conversations and suggestions.

This work was first presented by the second author at the Intensive Research Period on ``Euclidean Harmonic Analysis, Nilpotent Lie Groups and PDEs,'' Centro di Ricerca Matematica Ennio De Giorgi, Scuola Normale Superiore, Pisa, Italy, in April 2010.

\section{Notation and definitions}\label{s:notation}

Suppose that $\ap>0$ and $x$ is a point of $\R^n$, where $n\geq1$, $x=(x',x_n)$, $x'=(x_1,x_2,\ldots,x_{n-1})\in\R^{n-1}$ and $x_n\in\R$. Denote by $\Gamma_{\ap}(x)$ and $\tilde{\Gamma}_{\ap}(x)$ the cones of aperture $\ap$ and vertex $x$ given by
\[\Gamma_{\ap}(x)=\{y\in\R^n:\ap(y_n-x_n)>|y'-x'|\}\]
and
\[\tilde{\Gamma}_{\ap}(x)=\{y\in\R^n:\ap(x_n-y_n)>|y'-x'|\}.\]
The cone $\Gamma_{\ap}(x)$ lies above its vertex while $\tilde{\Gamma}_{\ap}(x)$ lies beneath.

A subset $\Omega$ of $\R^n$ is said to be a \textit{special Lipschitz domain} if 
\[\Omega=\{x\in\R^n:x_n>\lambda(x')\},\]
where the function $\lambda:\R^{n-1}\to\R$ satisfies the Lipschitz condition
\begin{equation}\label{Lipschitz}|\lambda(x')-\lambda(y')|\leq A|x'-y'|\qquad\forall x',y'\in\R^{n-1}\end{equation}
for some positive number $A$. The region that lies strictly beneath the graph of  $\lambda$ is denoted by $\Omega^-$. Thus $\Omega^-=\R^n\setminus\ov\Omega$. It follows immediately from the Lipschitz condition that when $0<\ap\leq A^{-1}$, the cone $\Gamma_{\ap}(x)$ is contained entirely in $\Omega$ whenever $x\in\overline{\Omega}$, while $\tilde{\Gamma}_{\ap}(x)$ lies entirely in $\Omega^-$ whenever $x\in\overline{\Omega^-}$.

Given a ball $B$ of $\R^n$, let $r(B)$ denote its radius and $|B|$ its volume. Whenever $c>0$, let $cB$ denote the ball with the same centre as $B$ and with radius $cr(B)$. If $r>0$ and $z\in \R^n$ then $B_r(z)$ denotes the ball $B$ of radius $r$ and centre $z$.

Throughout, let $\Cinf_0(\R^n)$ denote the space of smooth functions with compact support in $\R^n$. The Schwartz class of rapidly decreasing $\Cinf$ functions on $\R^n$ is denoted by $\mathcal{S}(\R^n)$. Given $f$ in $\mathcal{S}(\R^n)$, denote by $\widehat{f}$ its Fourier transform and by $\check{f}$ its inverse transform. Let $\mathcal{Z}(\R^n)$ denote the set of functions $\varphi$ in $\mathcal{S}(\R^n)$ possessing the property that $(\partial^{\alpha}\widehat{\varphi})(0)=0$ for every multi-index $\alpha$. Consider $\mathcal{Z}(\R^n)$ as a topological subspace of $\mathcal{S}(\R^n)$ and let $\mathcal{Z}'(\R^n)$ denote the topological dual of $\mathcal{Z}(\R^n)$. It is well-known (see, for example, \cite[Section 5.1.2]{hT83}) that $\mathcal{Z}'(\R^n)$ can be identified with the quotient space $\mathcal{S}'(\R^n)/\mathcal{P}(\R^n)$, where $\mathcal{S}'(\R^n)$ denotes the class of tempered distributions on $\R^n$ and $\mathcal{P}(\R^n)$ denotes the collection of all polynomials in $\R^n$. In this paper we use only the weak topology on $\mathcal{Z}'(\R^n)$. Hence we  say that a sequence $(u_n)_{n=1}^{\infty}\subset\mathcal{Z}'(\R^n)$ converges in $\mathcal{Z}'(\R^n)$ to an element $u$ of $\mathcal{Z}'(\R^n)$ if for each $\varphi$ in $\mathcal{Z}(\R^n)$, $|u_n(\varphi)-u(\varphi)|\to0$ as $n\to\infty$. The space $\mathcal{Z}'(\R^n)$ is complete with respect to this topology.

Many of the terms and definitions related to the space $\mathcal{S}'(\R^n)$ of tempered distributions have analogous formulations in $\mathcal{Z}'(\R^n)$. The Dirac delta distribution $\delta$ is defined as an element of $\mathcal{Z}'(\R^n)$ by the formula $\delta(\varphi)=\varphi(0)$ whenever $\varphi\in\mathcal{Z}(\R^n)$. Suppose that $u\in\mathcal{Z}'(\R^n)$. If $\Omega$ is an open set of $\R^n$ then we say that the support of $u$ is contained in $\ov\Omega$, or
$\supp u\subset\ov\Omega$, if $u(\varphi)=0$ for all $\phi$ in $\mathcal{Z}(\R^n)$ with support in $(\ov\Omega)^c$. 
 If $k\in\mathcal{S}(\R^n)$ then the convolution product $k\ast u\in\mathcal{Z}'(\R^n)$ is defined by
\[(k\ast u)(\varphi)=u(\tilde{k}\ast\varphi) \qquad\forall\varphi\in\mathcal{Z}(\R^n),\]
where $\tilde{k}(x)=k(-x)$ whenever $x\in\R$. (It is straightforward to show that $\tilde{k}\ast\varphi\in\mathcal{Z}(\R^n)$.) Moreover 
\begin{equation}
\supp (k\ast u) \subset \supp k + \supp u = \{x+y\in\R^n:\,x\in\supp k,\,y\in\supp u\}\ . \label{sppt}
\end{equation}

The homogeneous Besov spaces $\Bdot^s_{p,q}(\R^n)$ and homogeneous Triebel--Lizorkin spaces $\Fdot^s_{p,q}(\R^n)$ are subspaces of $\mathcal{Z}'(\R^n)$ defined in the following way. Fix a standard Littlewood--Paley dyadic system $(\theta_j)_{j\in\Z}$ of $\Cinf$ functions on $\R^n$ with support in dyadic annuli centred at the origin (see \cite[Chapter 5]{hT83} for details).  Given $u$ in $\mathcal{Z}'(\R^n)$, define $\Delta_j u$  by $\Delta_ju(\varphi)=u(\,(\theta_j\check{\varphi})\widehat{\,}\,)$ for all $\varphi$ in $\mathcal{Z}(\R^n)$. If $s\in\R$, $p>0$ and $q>0$ then  $\Bdot^s_{p,q}(\R^n)$ and $\Fdot^s_{p,q}(\R^n)$ are defined to be the spaces of all $u$ in $\mathcal{Z}'(\R^n)$ with finite norms given by 
\begin{align*}
\norm{u}_{\Bdot^s_{p,q}(\R^n)}&=\Big(\sum_{j\in\Z}2^{jsq}\norm{\Delta_j u}_{L^p(\R^n)}^q\Big)^{1/q}\qquad\text{and}\\
\norm{u}_{\Fdot^s_{p,q}(\R^n)}&=\Big\Vert\Big(\sum_{j\in\Z}2^{jsq}|\Delta_j u(\,\cdot\,)|^q\Big)^{1/q}\Big\Vert_{L^p(\R^n)}\ .
\end{align*}
By suitable modification, one may also define homogeneous spaces when $p=\infty$ or $q=\infty$.

 When $0<p<\infty$ and $0<q<\infty$, then $\mathcal{Z}(\R^n)$ is dense in  $\Bdot^s_{p,q}(\R^n)$ and  in $\Fdot^s_{p,q}(\R^n)$. 

Both these classes include the  homogeneous Sobolev spaces  $\Hdot^s(\R^n)$ as a special case, namely $\Hdot^s(\R^n)=\Bdot^s_{2,2}(\R^n)=\Fdot^s_{2,2}(\R^n)$ whenever $s\in\R$. Moreover the Hardy spaces $H^p(\R^n)$, which are defined as subspaces of $\mathcal{S}'(\R^n)$, can be characterised by such norms because the natural projection from $\mathcal{S}'(\R^n)$ to $\mathcal{Z}'(\R^n)$ induces an isomorphism from $H^p(\R^n)$ to $\Fdot^0_{p,2}(\R^n)$ when $0<p<\infty$ (see \cite[Section 5.2.4]{hT83}).

To simplify notation, given any real number $s$, let $\Adot^s(\R^n)$ denote any one of the spaces $\Bdot^s_{p,q}(\R^n)$ (for fixed $p$ and $q$ satisfying $0<p\leq\infty$, $0<q\leq\infty$) or $\Fdot^s_{p,q}(\R^n)$ (for fixed $p$ and $q$ satisfying $0<p<\infty$, $0<q\leq\infty$). Given a special Lipschitz domain $\Omega$ of $\R^n$, the spaces $\Adot^s_{\ov\Omega}(\R^n)$ and $\Adot^s(\Omega)$ are defined by
\begin{align*}
\Adot^s_{\ov\Omega}(\R^n)&=\{u\in\Adot^s(\R^n):\supp u\subset\ov\Omega\}\qquad\text{and}\\ \\
\Adot^s(\Omega)&=\Adot^s(\R^n)/\Adot^s_{\ov{\Omega^-}}(\R^n)\ .
\end{align*}
The quotient space $\Adot^s(\Omega)$ can be interpreted both as a space of `distributions' on $\R^n$ restricted to $\Omega$, and as a space of `distributions' on $\Omega$. The subspace $\Adot^s_{\ov\Omega}(\R^n)$ is a space of `distributions' on $\R^n$.

If $f\in\mathcal{S}(\R^n)$ and $t>0$ then $f_t$ is given by the formula
\[f_t(x)=t^{-n}f(x/t)\qquad\forall x\in\R^n.\]
It is easy to verify that the following formulae hold whenever $t>0$, $1\leq j\leq n$, and $f$ and $g$ belong to $\mathcal{S}(\R^n)$:
\[
(\partial_jf)_t=t\partial_j (f_t), \qquad
(f\ast g)_t=f_t\ast g_t, \quad\mbox{and}\quad
(f_t)\widehat{\,}(\xi)=\widehat{f}(t\xi).
\]

Denote the full exterior algebra on $\R^n$ by $\Wedge:=\Wedge^0\oplus\Wedge^1\oplus\dots \oplus\Wedge^n$, and the exterior product by $\wedg$. The interior product (or contraction) of $a = \sum_ja_jdx_j\in \Wedge^1$ with $u\in\Wedge^m$ is
$$
  a\contr u :
   = \sum_{k=1}^\ell 
     (-1)^{k-1} a_{j_k} dx_{j_1}\wedg\dots\wedg\widehat{dx}_{j_k}\wedg\dots\wedg dx_{j_\ell}\in\Wedge^{m-1}
$$
where the notation $\widehat{dx}_{j_k}$ means that the corresponding factor is to be omitted. (When $u\in\Wedge^0$, then $a\contr u =0$.)
The identity which we shall need is
\begin{equation} 
a\contr(b\wedg u)+b\wedg(a\contr u)=(a\cdot b)u\qquad\forall\, a,b\in\Wedge^1\ \text{and}\ u\in\Wedge\,. \label{eq:a contr b wedg w}
\end{equation}  

The space of differential forms with coefficients in $\Adot^s(\R^n)$ is denoted by $\Adot^s(\R^n,\Wedge)$, and similarly for forms whose coefficients lie in any of the test classes, functional classes, subspaces or quotient spaces mentioned above. The space of  forms in $\Adot^s(\R^n,\Wedge)$ of order $\ell$ is denoted by $\Adot^s(\R^n,\Wedge^{\ell})$. The topologies of these spaces are inherited in the obvious way. The exterior derivative $d$ for forms acting on $\R^n$ is defined by  $du=\Sigma_{j=1}^ne_j\wedg\partial_j u$, where $(e_j)_{j=1}^n$ denotes the standard basis for $\R^n$.  It has the property that $d^2=0$. By convention, $du=0$ whenever $u$ is an $n$-form. Note that $d$ maps $\mathcal{Z}(\R^n,\Wedge)$ into itself and hence is a well-defined map on $\mathcal{Z}'(\R^n,\Wedge)$. It also maps boundedly from $\Adot^s(\R^n,\Wedge)$ into $\Adot^{s-1}(\R^n,\Wedge)$.

\section{The potential map for special Lipschitz domains}\label{s:special lipschitz}

For the remainder of this paper, fix a Lipschitz constant $A$ and a constant $\sigma$ such that $A\geq0$, $\sigma>0$ and $\sigma A<1$. Suppose that $\Omega$ is any special Lipschitz domain, defined by a function $\lambda$ which satisfies the Lipschitz condition (\ref{Lipschitz}).   

 We also fix a function $\theta \in C^\infty(\R^n)$  with $\int_{\R^n} \theta(x)\,dx=1$ and 
\begin{equation}\label{set1}
\supp\theta \subset \{y\in\overline{\Gamma_\ap(0)}: \, 1\leq y_n, \, |y| \leq 2\}\ ,
\end{equation}
and define the vector valued $C^\infty$ function $\Theta$ with the same support by
\begin{equation}
\Theta(x) = \theta(x)\,x \qquad \forall x\in\R^n\ .
\end{equation} 
Before defining the potential map, we prove the following identities. 

\begin{proposition}\label{theta}
When $0<a<b<\infty$, define the $C^\infty$ function $\delta^{a,b}$ by
\begin{equation}
\delta^{a,b} := \int_a^b\,(\nabla\cdot \Theta)_t\,\frac{dt}t=\int_a^b\,\nabla\cdot \Theta_t\,dt\ .
\end{equation}
Then $\delta^{a,b} = \theta _a - \theta_b$.

For each $u\in{\mathcal Z}'(\R^n,\Wedge)$, $\delta^{a,b}\ast u \to u$ in ${\mathcal Z}'(\R^n,\Wedge)$ as $a\to0^+$ and $b\to\infty$. 
That is,
\[u=\int_0^\infty (\nabla\cdot\Theta)_t \ast u \frac{dt}{t}  \qquad \forall u\in\mathcal{Z}'(\R^n,\Wedge) \]
where the improper integral converges in $\mathcal{Z}'(\R^n,\Wedge)$.   
\end{proposition}

\begin{proof1} To verify  $\delta^{a,b} = \theta _a - \theta_b$, it suffices to check that, for $0<t<\infty$,
\begin{equation}\label{eq:rep}
\frac{d}{dt}\theta_t= -\frac1t(\nabla\cdot\Theta)_t\ .
\end{equation}
 A calculation using the product and chain rules shows that, for all $x\in\R^n$,
\begin{align}
\frac{d}{dt}\theta_t(x) =   \frac{d}{dt}\Big(t^{-n}\theta\big(\tfrac{x}{t}\big)\Big) 
&=t^{-n}\frac{d}{dt}\Big(\theta\big(\tfrac{x}{t}\big)\Big)-nt^{-n-1}\theta\big(\tfrac{x}{t}\big) \notag \\
&=-\frac{1}{t}\Big( \sum_{j=1}^n t^{-n}(\partial_j\theta)\big(\tfrac{x}{t}\big)\,\tfrac{x_j}{t} + nt^{-n}\theta(\tfrac{x}{t}\big) \Big) \notag \\
&=-\frac{1}{t}\Big( t^{-n} \sum_{j=1}^n (\partial_j\Theta_j)\big(\tfrac{x}{t}\big) \Big) \label{eq:rep calc}\notag \\
&=-\frac{1}{t}(\nabla\cdot\Theta)_t(x)\notag
\end{align}
as claimed. 

To complete the proof, we show that, for all $u\in {\mathcal Z}'(\R^n)$, (i) $\theta_t \ast u \to u$ as $t\to0^+$, and (ii) $\theta_t\ast u \to 0$ as $t\to \infty$.

To prove (i), it suffices to show that $\theta_t\ast\varphi\to\varphi$ in $\mathcal{Z}(\R^n)$ as $t\to0^+$ whenever $\varphi\in\mathcal{Z}(\R^n)$. Suppose that $\varepsilon>0$, $\varphi\in\mathcal{Z}(\R^n)$ and fix two multi-indices $\alpha$ and $\beta$. We need to show that 
\[\sup_{x\in\R^n}\left| x^\beta\,\partial^\alpha(\theta_t\ast\varphi-\varphi)(x) \right| <\varepsilon\]
whenever $t$ is sufficiently small. Since $\theta$ is compactly supported and $\varphi$ and all its derivatives have rapid decay, the supremum over the set $\{x\in\R^n:|x|>R\}$ can be made arbitrarily small by taking $R$ sufficiently large. Since $\theta$ has integral $1$, the family $(\theta_t)_{0<t<1}$ is an approximate identity, and hence the supremum over the ball $\{x\in\R^n:|x|\leq R\}$ can be made arbitrarily small  by taking $t$ sufficiently small. This completes the proof of (i).

We now prove (ii). Suppose that $\varepsilon>0$, $\varphi\in\mathcal{Z}(\R^n)$, and fix two multi-indices $\alpha$ and $\beta$. Since the Fourier transform is continuous in the topology of $\mathcal{S}(\R^n)$, its suffices to show that 
\[\sup_{\xi\in\R^n} \left| \xi^\beta \partial^\alpha(\widehat{\theta_t}\,\widehat{\varphi})(\xi) \right| < \varepsilon\]
whenever $t$ is sufficiently large. This is not difficult to achieve using the fact that $\widehat{\theta_t}(\xi)=\widehat\theta(t\xi)$, along with the assumption that $\widehat\varphi$ and 
all its derivatives are zero at $0\in\R^n$, as well as
the rapid decay of the functions $\widehat\theta$, $\widehat\varphi$ and their derivatives at $\infty$.
\end{proof1}

 \begin{remark}\label{convergence}
For those spaces $\Adot^{s}(\R^n)$ in which $\mathcal{Z}(\R^n)$ is dense,  the convergence $\delta^{a,b}\ast u \to u$ also holds in $\Adot^{s}(\R^n,\Wedge)$.  To see this, use the convergence result  in $\mathcal{Z}(\R^n,\Wedge)$, together with the uniform bound $\|\delta^{a,b}\ast u\|_{s}\leq C_s\|u\|_s$ for all $a<b$.
\end{remark} 

 We now define the {\it potential map} $T$. Whenever $0<a<b$, define the $\R^n$-valued kernel $K^{a,b}$ on $\R^n$ by
\begin{equation}\label{eq:truncated kernel1}
K^{a,b}(x)=\int_a^b \Theta_t(x) \,dt \qquad\forall x\in\R^n.
\end{equation}
For each $u$ in $\mathcal{Z}'(\R^n,\Wedge)$, define $Tu$ by the formula
\begin{equation}\label{eq:def of T}
Tu=\lim_{a\to0^+}\lim_{b\to\infty}K^{a,b}\ast\contr u,
\end{equation}
where both limits are taken in $\mathcal{Z}'(\R^n)$. If $u$ is a $0$-form then $Tu=0$. Note that  $T$ can be expressed by the formula
\[Tu(x)=\lim_{a\to0^+}\lim_{b\to\infty} \int_a^b \int_{\R^n} \theta_t(x-y)\,(x-y)\contr u(y)\,dy\,\frac{dt}{t}\]
whenever $u$ belongs to the subspace $\mathcal{Z}(\R^n,\Wedge)$ and where both limits are taken in $\mathcal{Z}'(\R^n)$.

We remark that in the case when $n=1$, the domain $\Omega$ is a semi-open interval $(\alpha, \infty)$, $A=0$ and $\Gamma_\sigma = (0,\infty)$. In this case, $\theta$ is a $C^\infty$ function with support in $[1, 2]$ such that $\int \theta(x)\, dx = 1$.
We then obtain that
\[Tu(x) = \int_{-\infty}^x e_1\contr u(y)\, dy\]
whenever $u\in\mathcal{Z}(\R,\Wedge)$ and $x\in\R$.  This is clearly the potential map for $d=\frac d{dx}\,e_1\wedg$ which preserves support in $(\alpha, \infty)$.

The following theorem, which is the main result of this section, asserts, among other things, that $T$ is well-defined. 

\begin{theorem}\label{th:T on special lipschitz}
Suppose that $\Omega$ is a special Lipschitz domain of $\R^n$ and that $s\in\R$. Then the operator $T$ defined above has the following properties:
\begin{enumerate}
\item[(i)] the operator $T$ maps from $\mathcal{Z}'(\R^n,\Wedge)$ to $\mathcal{Z}'(\R^n,\Wedge)$;
\item[(ii)]  the operator $T$ maps $\Adot^s(\R^n,\Wedge)$ to $\Adot^{s+1}(\R^n,\Wedge)$, and there is a constant $c$ such that
\[\norm{Tu}_{\Adot^{s+1}(\R^n,\Wedge)}\leq c\norm{u}_{\Adot^s(\R^n,\Wedge)} \qquad\forall u\in \Adot^s(\R^n,\Wedge);\]
\item[(iii)]  $dTu+Tdu=u$ whenever $u\in\mathcal{Z}'(\R^n,\Wedge)$;
\item[(iv)]  if  $u\in\mathcal{Z}'(\R^n,\Wedge)$ and $\supp(u)\subset\ov\Omega$ then $\supp(Tu)\subset\ov\Omega$; and 
\item[(v)] the operator $T$ maps $\Adot^s_{\ov\Omega}(\R^n,\Wedge)$ to $\Adot^{s+1}_{\ov\Omega}(\R^n,\Wedge)$, and for the same constant $c$ of Part (ii),
\[\norm{Tu}_{\Adot^{s+1}_{\ov\Omega}(\R^n,\Wedge)}\leq c\norm{u}_{\Adot^s_{\ov\Omega}(\R^n,\Wedge)} \qquad\forall u\in \Adot^s_{\ov\Omega}(\R^n,\Wedge)\,.\]
\end{enumerate}
\end{theorem}

Note that Part (v) is an immediate consequence of Parts (ii) and (iv) and the definition of $\Adot^s_{\ov\Omega}(\R^n,\Wedge)$. Before turning to the proof of the rest of the theorem, we give an immediate application to the regularity of the exterior derivative on special Lipschitz domains.

\begin{corollary}\label{cor:T on special lipschitz}
Suppose that $s\in\R$ and that $\Omega$ is a special Lipschitz domain of $\R^n$. If $u\in\Adot^s_{\ov\Omega}(\R^n,\Wedge)$ and $du=0$ then there exists $v$ in $\Adot^{s+1}_{\ov\Omega}(\R^n,\Wedge)$ such that $dv=u$. Moreover, there is a constant $c$ independent of $u$  such that
\[\norm{v}_{\Adot^{s+1}_{\ov\Omega}(\R^n,\Wedge)}\leq c \norm{u}_{\Adot^s_{\ov\Omega}(\R^n,\Wedge)}.\]
Consequently, the de Rham complex
\[ 0 \rightarrow \Adot^{s}_{\ov\Omega}(\R^n,\Wedge^0) \stackrel{d}{\rightarrow} \Adot^{s-1}_{\ov\Omega}(\R^n,\Wedge^1) \stackrel{d}{\rightarrow} \Adot^{s-2}_{\ov\Omega}(\R^n,\Wedge^2) \stackrel{d}{\rightarrow}\cdots \stackrel{d}{\rightarrow}\Adot^{s-n}_{\ov\Omega}(\R^n,\Wedge^n)  \rightarrow 0\]
is exact, and each space $\Adot^{s}_{\ov\Omega}(\R^n,\Wedge^k)$ has a direct sum decomposition
\[
\Adot^{s}_{\ov\Omega}(\R^n,\Wedge^k)=d\, \Adot^{s+1}_{\ov\Omega}(\R^n,\Wedge^{k-1}) \oplus Td\,\Adot^{s}_{\ov\Omega}(\R^n,\Wedge^k)
\]
with bounded projections $dT$ and $Td$.
\end{corollary}

\begin{proof1}
If $du=0$ then the anti-commutation relation in Part (iii) of the theorem becomes $dTu=u$. So set $v=Tu$. It is straightforward to check that $dT$ and $Td$ are bounded projections and that $dT\,\Adot^{s}_{\ov\Omega}(\R^n,\Wedge^k)=d\, \Adot^{s+1}_{\ov\Omega}(\R^n,\Wedge^{k-1})$. 
\end{proof1}

The important observation is that the projections are the same for all choices of the homogeneous Besov and Triebel--Lizorkin spaces $\Adot^{s}_{\ov\Omega}(\R^n,\Wedge^k)$, and consequently the spaces
$d\,\Adot^{s}_{\ov\Omega}(\R^n,\Wedge^k)$ of exact forms have the same interpolation properties as do the spaces $\Adot^{s}_{\ov\Omega}(\R^n,\Wedge^k)$.  

We remark that the operator $T$ as well as the constants in the estimates, depend on $A$ and the choice of $\ap$ and $\Theta$, but not on the precise domain $\Omega$. 

\begin{remark}\label{convergence2}
In the course of the proof, we show that the limit in formula \eqref{eq:def of T} holds in $\mathcal{Z}(\R^n,\Wedge)$ when $u\in \mathcal{Z}(\R^n,\Wedge)$, and in $\mathcal{Z}'(\R^n,\Wedge)$ when $u\in \mathcal{Z}'(\R^n,\Wedge)$. It can also be shown that, for those spaces $\Adot^{s}(\R^n)$ in which $\mathcal{Z}(\R^n)$ is dense,  the limit in formula \eqref{eq:truncated kernel1} holds in $\Adot^{s+1}(\R^n,\Wedge)$ when $u \in\Adot^{s}(\R^n,\Wedge)$. To see this, use the convergence result  in $\mathcal{Z}(\R^n,\Wedge)$, together with a uniform bound $\|K^{a,b}\ast u\|_{s+1}\leq C_s\|u\|_s$ for all $a<b$.
\end{remark} 

The rest of this section is devoted to proving Theorem~\ref{th:T on special lipschitz}. In Subsection~\ref{ss:homogeneous} we prove some elementary results about operators defined on $\mathcal{Z}'(\R^n)$ via homogeneous Fourier multipliers. These results are used in Subsection~\ref{ss:mapping properties of T} to prove Parts (i) and (ii) of Theorem~\ref{th:T on special lipschitz}. Parts (iii) and (iv) are then proved in Subsections~\ref{ss:anticommutation relations} and \ref{support properties}.

\subsection{Operators defined by homogeneous Fourier multipliers}\label{ss:homogeneous}

In this section we state and prove some rudimentary results about operators defined by homogeneous Fourier multipliers. These will be used later in the proof of Theorem \ref{th:T on special lipschitz}.

\begin{definition}
Suppose that  $k\in\R$. We say that a function $m$ on $\R^n$ is homogeneous of degree $k$ if
\begin{equation}\label{eq:homogeneous of degree k}
m(\tau\xi)=\tau^km(\xi) \qquad\forall\xi\in\R^n
\end{equation}
whenever $\tau>0$.
\end{definition}

\begin{proposition}\label{prop:hom}
Suppose that a function $m$ on $\R^n$ is homogeneous of degree $k$ and has partial derivatives of all orders on $\R^n\setminus\{0\}$. Then the operator $S$, given by
\[Su(\varphi)=u\big(\,(m\check{\varphi})\hat{\,}\,\big) \qquad\forall u\in\mathcal{Z}'(\R^n) \quad\forall \varphi\in\mathcal{Z}(\R^n)\ ,\]
is well-defined on $\mathcal{Z}'(\R^n)$, and maps each space $\Adot^s(\R^n)$ into $\Adot^{s-k}(\R^n)$ with
\[\norm{Su}_{\Adot^{s-k}(\R^n)} \leq c \norm{u}_{\Adot^s(\R^n)} \qquad\forall u\in\Adot^s(\R^n)\] 
for some constant $c=c(s)$.
\end{proposition}

\begin{proof1}
Assume the hypotheses of the proposition. It is obvious that, for each multi-index $\alpha$, $\partial^{\alpha}m$ is homogeneous of degree $k-|\alpha|$ and hence
\begin{equation}\label{eq:bound on partial m}
|\partial^{\alpha}m(\xi)|\leq c|\xi|^{-|\alpha|+k} \qquad\forall\xi\in\R^n\setminus\{0\}
\end{equation}
where $c=\sup\{|\partial^{\alpha}m(\omega)|:|\omega|=1\}$. 

To show that $S$ is well-defined on $\mathcal{Z}'(\R^n)$, it suffices to verify that $(m\check{\varphi})\,\hat{\,}\in\mathcal{Z}(\R^n)$ whenever $\varphi\in\mathcal{Z}(\R^n)$. Since $m$ has at most polynomial growth at infinity, one need only show that $\partial^{\alpha}(m\check{\varphi})(0)=0$ for every multi-index $\alpha$. But this follows from (\ref{eq:bound on partial m}) and the fact that $\check{\varphi}(\xi)$ is $O(|\xi|^N)$ for every positive integer $N$.

To prove the bound, consider the Fourier multiplier $m_k$ given by $m_k(\xi)=|\xi|^{-k}m(\xi)$. Clearly $m_k$ has derivatives of all orders away from $0$ and is homogeneous of degree $0$. So the operator $Q$, defined by
\[Qu(\varphi)=u\big(\,(m_k\check{\varphi})\hat{\,}\,\big) \qquad\forall u\in\mathcal{Z}'(\R^n) \quad\forall\varphi\in\mathcal{Z}(\R^n)\]
is well-defined on $\mathcal{Z}'(\R^n)$. Since
\[\sup\left\{|\xi|^{|\alpha|}|\partial^{\alpha}m_k(\xi)|:\xi\in\R^n\setminus\{0\},|\alpha|\leq N\right\}\]
is bounded for each positive integer $N$, the operator $Q$ is bounded on $\Adot^s(\R^n)$ by standard Fourier multiplier theory (see, for example, \cite[Theorem 5.2.2]{hT83}). Now $S=\dot{I}_{k}Q$, where the lifting operator $\dot{I}_{k}$, given by
\[\dot{I}_{k}u(\varphi)=u\big(\,(|\iota|^k\check{\varphi})\hat{\,}\,\big) \qquad\forall u\in\mathcal{Z}'(\R^n)\quad\forall \varphi\in\mathcal{Z}(\R^n),\]
maps $\Adot^s(\R^n)$ isomorphically onto $\Adot^{s-k}(\R^n)$. Here $\iota$ denotes the $\R^n$-valued function on $\R^n$ given by $\iota(\xi)=\xi$.
This completes the proof.
\end{proof1}

\subsection{Mapping properties of $T$}\label{ss:mapping properties of T}

Our aim in this subsection is twofold: first, to show that the operator $T$ defined by the limit (\ref{eq:def of T}) is well-defined; and second, to prove Parts (i) and (ii) of Theorem~\ref{th:T on special lipschitz}.

In fact, we prove these results within a more general setting. Given a function $\psi$ in $\Cinf_0(\R^n)$, define the truncated kernel $k^{a,b}$ by
\[k^{a,b}(x)=\int_a^b\psi_t(x)\,dt \qquad\forall x\in\R^n\]
whenever $0<a<b$. Since $k^{a,b}\in\Cinf_0(\R^n)$, the operator $S^{a,b}$, given by
\[S^{a,b}u=k^{a,b}\ast u \qquad\forall u\in\mathcal{Z}'(\R^n)\]
is well-defined on $\mathcal{Z}'(\R^n)$. Denote $\widehat{k^{a,b}}$ by $m^{a,b}$.

Note that each component of $K^{a,b}*\contr u$, where $K^{a,b}$ is the kernel used to define $T$, is of the form $S^{a,b}u_I$, where $\psi(x)=\theta(x)x_j$ and $u_I$ is a component of $u$ corresponding to an index set $I$ including $j$.

The first lemma of this subsection will be used to show that the limit
\[\lim_{a\to0^+}\lim_{b\to\infty}S^{a,b}u,\]
taken in the topology of $\mathcal{Z}'(\R^n)$, is well-defined whenever $u\in\mathcal{Z}'(\R^n)$.

\begin{lemma}\label{lem:uniform convergence of multipliers} Suppose that $(m^{a,b})_{0<a<b}$ is the net of Schwartz functions defined above. Then the function $m$, given by $m(0)=0$ and
\begin{equation}\label{eq:m}
m(\xi)=\lim_{a\to0^+}\lim_{b\to\infty}m^{a,b}(\xi) \qquad\forall\xi\in\R^n\setminus\{0\},
\end{equation}
is well-defined. Moreover, $m\in\Cinf(\R^n\setminus\{0\})$ and for any multi-index $\alpha$,
\begin{equation}\label{eq:partial m}
\partial^{\alpha}m(\xi)=\lim_{a\to0^+}\lim_{b\to\infty}\partial^{\alpha}m^{a,b}(\xi) \qquad\forall\xi\in\R^n\setminus\{0\},
\end{equation}
where the convergence is uniform on annuli centred at the origin. Finally, $m$ is homogeneous of degree $-1$ on $\R^n$.
\end{lemma}

\begin{proof1}
Fix a multi-index $\alpha$. To show that $m$ is well-defined and its derivatives are given by (\ref{eq:partial m}), it suffices to show that the net $(\partial^{\alpha}m^{a,b})_{a,b}$ is uniformly Cauchy on the annulus $\{\xi\in\R^n:r\leq|\xi|\leq R\}$ for some fixed numbers $r$ and $R$ satisfying $0<r<R<\infty$. Henceforth, suppose that $r\leq|\xi|\leq R$.

To begin, note that
\[m^{a,b}(\xi)=\int_a^b\widehat{\psi}(t\xi)\,dt\]
whenever $0<a<b$. Therefore
\begin{align*}
\left|(\partial^{\alpha}m^{a,b})(\xi)\right| 
&=\left| \int_{a}^{b}t^{|\alpha|}(\partial^{\alpha}\widehat{\psi})(t\xi)\,dt\right| \\
&\leq |\xi|^{-|\alpha|}\int_{a}^{b}|t\xi|^{|\alpha|}|(\partial^{\alpha}\widehat{\psi})(t\xi)|\,dt \\
&\leq c_{\alpha}|\xi|^{-|\alpha|}\int_{a}^{b}\min\big\{1,|t\xi|^{-2}\big\}\,dt 
\end{align*}
for some constant $c_{\alpha}$, since $\widehat{\psi}\in\mathcal{S}(\R^n)$. Now
\[\left|(\partial^{\alpha}m^{b,b_0})(\xi)\right|  \leq c_{\alpha} |\xi|^{-|\alpha|-2}\left(\frac{1}{b}-\frac{1}{b_0}\right)<\frac{c_{\alpha}}{r^{|\alpha|+2}b}\]
whenever $b_0>b$, while
\[\left|(\partial^{\alpha}m^{a_0,a})(\xi)\right|  \leq c_{\alpha} |\xi|^{-|\alpha|}(a-a_0)<c_{\alpha}r^{-|\alpha|}a\]
whenever $a_0<a$. So if $0<a_0<a<b<b_0$ then 
\[\left|(\partial^{\alpha}m^{a_0,b_0})(\xi)-(\partial^{\alpha}m^{a,b})(\xi)\right|<c_{\alpha}\left(r^{-|\alpha|-2}b^{-1}+r^{-|\alpha|}a\right),\]
which can be made as small as we like by taking $b$ sufficiently large and $a$ sufficiently close to $0$. 

The fact that $m$ is homogeneous of degree $-1$ follows easily from the homogeneity of each $m^{a,b}$ and the definition of $m$. This completes the proof.
\end{proof1}

Using Proposition \ref{prop:hom} and Lemma~\ref{lem:uniform convergence of multipliers}, we now define the operator $S$ on $\mathcal{Z}'(\R^n)$ by
\[Su(\varphi)=u\big(\,(m\check{\varphi})\hat{\,}\,\big) \qquad\forall u\in\mathcal{Z}'(\R^n)\quad\forall\varphi\in\mathcal{Z}(\R^n),\]
where $m$ is the function given by (\ref{eq:m}).

\begin{lemma}
If $u\in\mathcal{Z}'(\R^n)$ then $S^{a,b}u$ converges to $Su$ in  $\mathcal{Z}'(\R^n)$ as $a\to0^+$ and as $b\to\infty$.
\end{lemma}

\begin{proof1}
Suppose that $\varphi\in\mathcal{Z}(\R^n)$ and $u\in\mathcal{Z}'(\R^n)$. Since
\[S^{a,b}u(\varphi)=u\big(\,(m^{a,b}\check{\varphi})\hat{\,}\,\big) \qquad\forall u\in\mathcal{Z}'(\R^n)\quad\forall\varphi\in\mathcal{Z}(\R^n),\]
it suffices to show that $m^{a,b}\check{\varphi}\to m\check{\varphi}$ in $\mathcal{S}(\R^n)$ as $a\to0^+$ and $b\to\infty$. For then $(m^{a,b}\check{\varphi})\hat{\,} \to (m\check{\varphi})\hat{\,}$ in $\mathcal{S}(\R^n)$ and hence in $\mathcal{Z}(\R^n)$. It follows that
\[(S^{a,b}u-Su)(\varphi)=u\left(\big((m^{a,b}-m)\check{\varphi}\big)\hat{\phantom{l}}\,\right)\to0\]
as $a\to0^+$ and $b\to\infty$, which establishes the lemma.

Suppose that $\alpha$ and $\beta$ are two multi-indices and $\epsilon>0$. We need to show that there exist positive numbers $a_0$ and  $b_0$ such that
\[
\sup_{\xi\in\R^n}\left|\xi^{\alpha}\partial^{\beta}\big((m-m^{a,b})\check{\varphi}\big)(\xi)\right|<\epsilon
\]
whenever $0<a<a_0<b_0<b$. By expanding the left-hand side using the multidimensional version of Leibniz' rule, it suffices to show that there are positive numbers $a_0$ and  $b_0$ such that
\begin{equation}\label{eq:seminorm estimate}
\sup_{\xi\in\R^n}\,\,\sup_{|\gamma|\leq|\beta|}\, |\xi|^{|\alpha|}\big|\partial^{\gamma}(m-m^{a,b})(\xi)\big|\,\big|\partial^{\beta-\gamma}\check{\varphi}(\xi)\big|<\frac{\epsilon}{c_{\beta}}
\end{equation}
whenever $0<a<a_0<b_0<b$, where the constant $c_{\beta}$ is the largest coefficient appearing the formula for Leibniz' rule. 

By Lemma~\ref{lem:uniform convergence of multipliers}, there exist $a_1$ and $b_1$ such that
\begin{equation}\label{eq:c difference}
\sup\left\{|\partial^{\gamma}(m-m^{a,b})(\omega)|:\,|\gamma|\leq|\beta|,\,|\omega|=1\right\}\leq1
\end{equation}
whenever $0<a<a_1<b_1<b$. Since $m$ and each $m^{a,b}$ are homogeneous of degree $-1$, $m-m^{a,b}$ is also homogeneous of degree $-1$ and consequently (\ref{eq:bound on partial m}) and (\ref{eq:c difference}) give the estimate
\begin{equation}\label{eq:m difference}
|\partial^{\gamma}(m-m^{a,b})(\xi)|\leq |\xi|^{-|\gamma|-1} \qquad\forall\xi\in\R^n\setminus\{0\}
\end{equation}
whenever $0<a<a_1'<b_1<b$ and $|\gamma|\leq|\beta|$. Now choose $R$ in $(0,\infty)$ so large that
\begin{equation}\label{eq:m R}
\sup_{|\xi|>R}\,\,\sup_{|\gamma|\leq|\beta|}\,|\xi|^{|\alpha|-|\gamma|-1}|\partial^{\beta-\gamma}\check{\varphi}(\xi)|<\frac{\epsilon}{3c_{\beta}}\,.
\end{equation}
This is possible since $\partial^{\beta-\gamma}\check{\varphi}$ is rapidly decreasing at infinity. Now choose $r$ in $(0,\infty)$ so small that
\begin{equation}\label{eq:m r}
\sup_{0<|\xi|<r}\,\,\sup_{|\gamma|\leq|\beta|}\,\xi|^{|\alpha|-|\gamma|-1}|\partial^{\beta-\gamma}\check{\varphi}(\xi)|<\frac{\epsilon}{3c_{\beta}}\,.
\end{equation}
This is possible since $\check{\varphi}$ and all its partial derivatives are $0$ at the origin. By Lemma~\ref{lem:uniform convergence of multipliers}, there are positive numbers $a_2$ and $b_2$ such that
\begin{equation}\label{eq:m r R}
\sup_{r\leq|\xi|\leq R}\,\,\sup_{|\gamma|\leq|\beta|}\,\left|\xi^{\alpha}\,\big(\partial^{\gamma}(m-m^{a,b})\big)(\xi)\,\big(\partial^{\beta-\gamma}\check{\varphi}\big)(\xi)\right|<\frac{\epsilon}{3c_{\beta}}
\end{equation}
whenever $0<a<a_2<b_2<b$. By combining estimates (\ref{eq:m difference}), (\ref{eq:m R}), (\ref{eq:m r}) and (\ref{eq:m r R}), it is easy to see that  (\ref{eq:seminorm estimate}) holds whenever $0<a<\min(a_0,a_1)<\max(b_0,b_1)<b$. This shows that $m^{a,b}\check{\varphi}\to m\check{\varphi}$ in $\mathcal{S}(\R^n)$ and completes the proof.
\end{proof1}

We observe now that we have proved Parts (i) and (ii) of Theorem~\ref{th:T on special lipschitz}. Each component of $T$ is a limit in $\mathcal{Z}'(\R^n)$ of convolution operators of the form $S^{a,b}$ and therefore has the same properties as $S$. In particular, Proposition~\ref{prop:hom} shows that $T$ maps from $\mathcal{Z}'(\R^n,\Wedge)$ into $\mathcal{Z}'(\R^n,\Wedge)$ and boundedly lifts `functions' of degree $s$ in the homogeneous Besov and Triebel--Lizorkin spaces to `functions' of degree $s+1$.

\subsection{Anticommutation relations}\label{ss:anticommutation relations}

We now turn to the proof of Theorem~\ref{th:T on special lipschitz} (iii). Suppose that $0<a<b$. Define the operator $T^{a,b}$ on $\mathcal{Z}'(\R^n,\Wedge)$ by $T^{a,b}u=K^{a,b}\ast\contr u$ whenever $u\in\mathcal{Z}'(\R^n,\Wedge)$, where $K^{a,b}$ denotes the $\R^n$-valued kernel given by (\ref{eq:truncated kernel1}). Recall that $\delta^{a,b}$ denotes the function in $\mathcal{Z}(\R^n)$ given by
\[\delta^{a,b} = \int_a^b\,(\nabla\cdot \Theta)_t\,\frac{dt}t = \int_a^b\,\nabla\cdot \Theta_t\,dt\ .
\]
Fix $u$ in $\mathcal{Z}'(\R^n,\Wedge)$. If we can show that
\begin{equation}\label{eq:anticommutation a,b}
dT^{a,b}u+T^{a,b}du=\delta^{a,b}\ast u\ ,
\end{equation}
then, by taking limits as $a\to0^+$ and $b\to\infty$ in $\mathcal{Z}'(\R^n)$ and applying Proposition \ref{theta}, Part (iii) of Theorem~\ref{th:T on special lipschitz} will be proved.

Now 
\begin{align*}
T^{a,b}u&=\int_a^b\Theta_t\ast \contr u\,dt\,,
\intertext{so}
dT^{a,b}u&=\sum_{j=1}^n e_j\wedg\partial_j\int_a^b\Theta_t\ast \contr u\,dt\\
&=\sum_{j=1}^n \int_a^b e_j\wedg(\partial_j\Theta_t\ast \contr u)\,dt
\intertext{and}
T^{a,b}du&=\int_a^b\Theta_t\ast \contr du\,dt\\
&=\sum_{j=1}^n \int_a^b\partial_j  \Theta_t\ast \contr(e_j\wedg u)\,dt\, .
\end{align*}
Therefore, using the identity (\ref{eq:a contr b wedg w}), we obtain
\begin{align*}
dT^{a,b}u+T^{a,b}du&=\sum_{j=1}^n \int_a^b\partial_j \Theta_t \cdot e_j\ast u\,dt\\
&=\int_a^b\nabla \cdot \Theta_t\ast u\,dt\\
&=\delta^{a,b}\ast u
\end{align*}
as required.

Hence we have shown equation (\ref{eq:anticommutation a,b}). This completes the proof of Theorem~\ref{th:T on special lipschitz} (iii).

\subsection{Support properties of $T$}\label{support properties}

To complete the proof of Theorem~\ref{th:T on special lipschitz}, it remains to show Part (iv). 
Suppose $\supp u \subset \overline\Omega$. By (\ref{sppt}) and (\ref{set1}),
\begin{align*}
\supp(\Theta_t\ast\contr u)
&\subset  \supp\Theta_t+\supp u\\
&\subset\overline{\Gamma_\ap(0)} +\overline{\Omega} \ \subset\ \overline\Omega
\end{align*}
whenever $t\in(0,\infty)$. Hence  $\supp T^{a,b}u\subset \ov\Omega$ whenever $0<a<b$. By taking limits in $\mathcal{Z}'(\R^n)$ as $a\to0^+$ and $b\to\infty$, one concludes that $\supp(Tu)\subset\ov\Omega$.

This completes the proof of Theorem~\ref{th:T on special lipschitz}.{\hfill \ensuremath{\Box}}

\section{Analogous results for complementary domains and restriction spaces}\label{s:restriction spaces}

Suppose that $\Omega$ is a special Lipschitz domain of $\R^n$, and recall that $\Omega^-$ denotes the region strictly below the corresponding Lipschitz graph. Define an operator $\tilde{T}$ by the formula
\[\tilde{T}u=\lim_{a\to0^+}\,\lim_{b\to\infty}\tilde{K}^{a,b}\ast\contr\, u, \qquad\forall u\in\mathcal{Z}'(\R^n,\Wedge),\]
where both limits are taken in $\mathcal{Z}'(\R^n)$ and $\tilde{K}^{a,b}(x)=K^{a,b}(-x)$ for all $x$ in $\R^n$. Here $K^{a,b}$ is the truncated kernel given by (\ref{eq:truncated kernel1}). The analytic properties of $\tilde{T}$ are clearly the same as those of $T$. However, if $\supp u\subset\ov{\Omega^-}$ then $\supp(\tilde{T}u)\subset\ov{\Omega^-}$. Hence, in Theorem \ref{th:T on special lipschitz}, one may replace $T$ by $\tilde{T}$ throughout, and $\Omega$ by $\Omega^-$ in Parts (iv) and (v), to obtain an analogous result for the complementary Lipschitz domain $\Omega^-$.

We now draw some conclusions for the restriction space $\Adot^s(\Omega,\Wedge)$. Given $u$ in $\Adot^s(\R^n,\Wedge)$, let $[u]$ denote the equivalence class with representative $u$ associated to the equivalence relation
\[v\sim w \iff v-w\in\Adot^s_{\ov{\Omega^-}}(\R^n,\Wedge).\]
By definition, $[u]$ belongs to $\Adot^s(\Omega,\Wedge)$ and conversely every element of $\Adot^s(\Omega,\Wedge)$ is of this form. Define an operator $R$ by
\[R[u]=[\tilde{T}u] \qquad \forall \,[u]\in\Adot^s(\Omega,\Wedge).\]
Since $\tilde{T}$ maps boundedly from $\Adot^s(\R^n,\Wedge)$ to $\Adot^{s+1}(\R^n,\Wedge)$ and from $\Adot^s_{\ov{\Omega^-}}(\R^n,\Wedge)$ to $\Adot^{s+1}_{\ov{\Omega^-}}(\R^n,\Wedge)$, the operator $R$ is well-defined and maps boundedly from $\Adot^s(\Omega,\Wedge)$ to $\Adot^{s+1}(\Omega,\Wedge)$.

Similarly, the exterior derivative $d$ is defined as an operator on $\Adot^s(\Omega,\Wedge)$ by
\[d[u]=[du] \qquad \forall \,[u]\in\Adot^s(\Omega,\Wedge),\]
and maps boundedly from $\Adot^s(\Omega,\Wedge)$ into $\Adot^{s-1}(\Omega,\Wedge)$. We thus obtain another variant of Theorem~\ref{th:T on special lipschitz}.

\begin{proposition}\label{prop:R}
Suppose that $\Omega$ is a special Lipschitz domain of $\R^n$ and that $s\in\R$. Then the operator $R$ defined above has the following properties:
\begin{enumerate}
\item[(i)]  the operator $R$ maps $\Adot^s(\Omega,\Wedge)$ to $\Adot^{s+1}(\Omega,\Wedge)$, and there is a constant $c$ such that
\[\norm{R[u]}_{\Adot^{s+1}(\Omega,\Wedge)}\leq c\norm{[u]}_{\Adot^s(\Omega,\Wedge)} \qquad\forall [u]\in \Adot^s(\Omega,\Wedge);\]
\item[(ii)]  $dR[u]+Rd[u]=[u]$ whenever $[u]\in\Adot^s(\Omega,\Wedge)$.
\end{enumerate}
\end{proposition}

One immediately obtains a regularity result for the exterior derivative on $\Adot^s(\Omega,\Wedge)$.

\begin{corollary}\label{cor:R1}
Suppose that $s\in\R$ and $\Omega$ is a special Lipschitz domain. If $[u]\in\Adot^s(\Omega,\Wedge)$ and $d[u]=0$ then there exists $[v]$ in $\Adot^{s+1}(\Omega,\Wedge)$ and a constant $c$ independent of $[u]$ such that $d[v]=[u]$ and
\[\norm{[v]}_{\Adot^{s+1}(\Omega,\Wedge)} \leq c \norm{[u]}_{\Adot^s(\Omega,\Wedge)}.\]
Consequently, the de Rham complex
\[ 0 \rightarrow \Adot^{s}(\Omega,\Wedge^0) \stackrel{d}{\rightarrow} \Adot^{s-1}(\Omega,\Wedge^1) \stackrel{d}{\rightarrow} \Adot^{s-2}(\Omega,\Wedge^2) \stackrel{d}{\rightarrow}\cdots \stackrel{d}{\rightarrow}\Adot^{s-n}(\Omega,\Wedge^n)  \rightarrow 0\]
is exact, and each space $\Adot^{s}(\Omega,\Wedge^k)$ has a direct sum decomposition
\[\Adot^{s}(\Omega,\Wedge^k)=d\, \Adot^{s+1}(\Omega,\Wedge^{k-1}) \oplus Rd\,\Adot^{s}(\Omega,\Wedge^k)\]
with bounded projections $dR$ and $Rd$.
\end{corollary}

\section{Reproducing formula of Chang, Krantz and Stein}\label{sec:CKS}

 Our construction of the potential map  $T$ was motivated by the reproducing formula used by Chang, Krantz and Stein in obtaining an atomic decomposition of functions in a Hardy space on a special Lipschitz domain. We indicate here the connection between our Proposition \ref{theta} and their
  reproducing formula \cite[Lemma 3.4, Lemma 3.5]{CKS93}.   

\begin{proposition}\label{phi}
Let $\phi\in C^\infty(R^n)$  with $\int_{\R^n} \phi(x)\,dx=1$ and 
\begin{equation}\label{set}
\supp\phi \subset \{y\in\overline{\Gamma_\ap(0)}: \, \tfrac12\leq y_n, \, |y| \leq 1\}\ ,
\end{equation}
and define the vector valued $C^\infty$ function $\Psi$ with the same support by
\begin{equation}
\Psi(x) = 2\phi(x)\,x \qquad \forall x\in\R^n\ .
\end{equation}
When $0<a<b<\infty$, define the $C^\infty$ function $\delta^{a,b}$ by
\begin{equation}
\delta^{a,b} := \sum_{j=1}^n\int_a^b\,(\partial_j\phi)_t\ast(\Psi_j)_t\,\frac{dt}t\ .\end{equation}
Then, for each $u\in{\mathcal Z}'(\R^n,\Wedge)$, $\delta^{a,b}\ast u \to u$ in ${\mathcal Z}'(\R^n,\Wedge)$ as $a\to0^+$ and $b\to\infty$.
 \end{proposition}

\begin{proof1}
The function $\theta: = \phi\ast\phi$ satisfies the hypotheses of Proposition \ref{theta}, so
$\delta^{a,b}\ast u \to u$ in ${\mathcal Z}'$ for all $u\in{\mathcal Z}'$ where 
\begin{equation}
\delta^{a,b} := \int_a^b\,(\nabla\cdot \Theta)_t\,\frac{dt}t
\end{equation}
and $\Theta(x):=\theta(x)x$. We shall show that $\phi\ast\Psi=\Theta$, and hence $\sum_{j=1}^n \partial_j\phi\ast\Psi_j=\nabla\cdot \Theta$, from which it follows by scaling and integrating that 
the two expressions  for $\delta^{a,b}$ are the same, thus proving the result.

What needs to be shown is that $(\phi\ast\Psi)(x)=(\phi\ast\phi)(x)x$ for all $x\in \R^n$. Indeed,
\begin{align*}
(\phi\ast\Psi)(x)
&=\int_{R^n}\phi(x-y)2\phi(y)y\,dy\\
&=\int_{\R^n}\phi(x-y)\phi(y)y\,dy+\int_{\R^n}\phi(y)\phi(x-y)(x-y)\,dy\\
&=(\phi\ast\phi)(x) x\ ,
\end{align*}
thus proving the proposition.
\end{proof1}

We remark that, to be useful as a reproducing formula, the function $\phi$ can be chosen with some zero  moments, in particular \begin{equation}\label{meanzero}\int_{\R^n}\Psi(x)\,dx  = 2 \int_{\R^n}\phi(x)\,x\,dx  = 0
\ .\end{equation} 

With the above choice of functions, the operator $T$ has the form
\begin{align}\label{reproducing}
Tu&=\lim_{a\to0^+}\lim_{b\to\infty}\int_a^b\Theta_t\ast\contr u\,dt\notag\\
&=\lim_{a\to0^+}\lim_{b\to\infty}\int_a^b\phi_t\ast\Psi_t\ast\contr u\,dt\ .
\end{align}

As before, convolution with $\phi_t$ and $\Psi_t$ preserve support in the special Lipschitz domain $\Omega$, and $T$ is a bounded operator from ${\Adot}^s(\overline\Omega,\Wedge)$ to ${\Adot}^{s+1}(\overline\Omega,\Wedge)$ for every choice of ${\Adot}^s$. Moreover, by Remark \ref{convergence2}, for those spaces ${\Adot}^{s}(\R^n)$ in which ${\mathcal Z}(\R^n)$ is dense, the limits exist in ${\Adot}^{s+1}(\overline\Omega,\Wedge)$ whenever $u\in {\Adot}^{s}(\overline\Omega,\Wedge)$.

\section{Atomic decomposition of Hardy spaces of exact forms on special Lipschitz domains}\label{s:atomic decomp}

In this section we use the operator $T$ of Theorem~\ref{th:T on special lipschitz} and the reproducing formulae  above, including the zero moment condition \eqref{meanzero}, to show that Hardy spaces of exact forms on special Lipschitz domains can be characterised by atomic decompositions. 

In the following definitions of these spaces  and their corresponding atoms, we at first allow $\Omega$ to be an arbitrary domain in $\R^n$, where $n\geq1$. 

\begin{definition}\label{def:H^p_d}
Suppose that $1\leq\ell\leq n$ and $n/(n+1)<p\leq1$. Let $H^p_d(\R^n,\Wedge^{\ell})$ denote the Hardy  space of all $\ell$-forms $u$ in $H^p(\R^n,\Wedge)$ such that $u=dv$ for some $(\ell-1)$-form $v$ in $\mathcal{S}'(\R^n,\Wedge^{\ell-1})$. Given a domain $\Omega$ in $\R^n$, we say that $u$ is in $H^p_{z,d}(\ov\Omega,\Wedge^{\ell})$ if $u\in H^p_d(\R^n,\Wedge^{\ell})$ and there exists $v$ in $\mathcal{S}'(\R^n,\Wedge^{\ell-1})$ such that $u=dv$ and $\supp v\subset\ov\Omega$.
\end{definition}

\begin{remark}\label{rem:H1}
Definition~\ref{def:H^p_d} was first introduced in the papers \cite{LM05} and \cite{LM04} of Lou and McIntosh for the case when $p=1$. When $n/(n+1)<p\leq1$, the space $H^p_d(\R^n,\Wedge^n)$ is isomorphic to the classical Hardy space $H^p(\R^n)$, while  $H^p_{z,d}(\ov\Omega,\Wedge^n)$ is isomorphic to the Hardy space $H^p_z(\Omega)$ of Chang, Krantz and Stein \cite{CKS93}. 
\end{remark}

Following \cite{LM05} and \cite{LM04}, we introduce atoms of Hardy spaces of exact forms.

\begin{definition}\label{def:H^1 atom}
Suppose that $1\leq\ell\leq n$ and $n/(n+1)<p\leq1$. We say that $a$ is an $H^p_d(\R^n,\Wedge^{\ell})$-atom if for some ball $B$ in $\R^n$,
\begin{enumerate}
\item[(a)]  there exists $b$ in $L^2(\R^n,\Wedge^{\ell-1})$ such that $\supp b\subset \ov B$ and $a=db$, and
\item[(b)] $\ds\norm{a}_{L^2(\R^n,\Wedge)}\leq|B|^{1/2-1/p}$.
\end{enumerate}
\end{definition}

Note that if $n/(n+1)<p\leq 1$ and $a$ is an $H^p_d(\R^n,\Wedge^{\ell})$-atom, then each component of $a$ is a classical $H^p(\R^n)$-atom.

\begin{definition}\label{def:H^1_omega atom}
Suppose that $\Omega$ is a domain of $\R^n$, $1\leq\ell\leq n$ and $n/(n+1)<p\leq 1$. We say that $a$ is an $H^p_{z,d}(\ov\Omega,\Wedge^{\ell})$-atom if for some ball $B$ in $\R^n$,
\begin{enumerate}
\item[(a)]  there exists $b$ in $L^2(\R^n,\Wedge^{\ell-1})$ such that $\supp b\subset\ov B$ and $a=db$, 
\item[(b)] $\ds\norm{a}_{L^2(\R^n,\Wedge)}\leq|B|^{1/2-1/p}$, and
\item[(c)] $4B\subset\Omega$.
\end{enumerate}
\end{definition}

Note that, following \cite{LM04}, the supports of $H^p_{z,d}(\ov\Omega,\Wedge^{\ell})$-atoms are away from the boundary of $\Omega$, which is stronger than the classical definition of \cite{CKS93}. 

The following lemma gives an $L^2$ estimate for the function $b$ of Definitions~\ref{def:H^1 atom} and \ref{def:H^1_omega atom}.

\begin{lemma}\label{lem:b bound}
Suppose that $1\leq\ell\leq n$,  $n/(n+1)<p\leq 1$ and $a$ is an $H^p_d(\R^n,\Wedge^{\ell})$-atom (respectively an $H^p_{z,d}(\ov\Omega,\Wedge^{\ell})$-atom). Then the $L^2(\R^n,\Wedge^{\ell-1})$ form $b$ of Definition~\ref{def:H^1 atom} (respectively Definition~\ref{def:H^1_omega atom}) can be chosen such that
\[\norm{b}_{L^2(\R^n,\Wedge)}\leq c_n r(B) |B|^{1/2-1/p}\]
where the constant $c_n$ depends only on $n$.
\end{lemma}

\begin{proof1}
Let $\mathscr{B}$ denote the collection of all balls in $\R^n$ and let $L^2_{\ov B}(\R^n,\Wedge^k)$ denote the space of $k$-forms with components in $L^2(\R^n)$ and support in the closure of a ball $B$. By applying the result of \cite{mM04} or \cite[Section 3]{CoMcI} to a unit ball and then scaling, one obtains the following. There exists a constant $c_n$ and a family of operators $\{T_{\ell}^B:1\leq \ell\leq n, B\in\mathscr{B}\}$ with the following properties:
\begin{itemize}
\item[(i)] if $B\in\mathscr{B}$ and $1\leq\ell\leq n$ then $T^B_{\ell}$ maps from $L^2_{\ov B}(\R^n,\Wedge^{\ell})$ to $L^2_{\ov B}(\R^n,\Wedge^{\ell-1})$ and
\[\norm{T_{\ell}^Bu}_{L^2(\R^n,\Wedge)} \leq c_nr(B) \norm{u}_{L^2(\R^n,\Wedge)} \qquad\forall u\in L^2_{\ov B}(\R^n,\Wedge^{\ell});\]
\item[(ii)] if $B\in\mathscr{B}$ and $1\leq\ell< n$ then $dT^B_{\ell}u+T^B_{\ell+1}du=u$ for every $u\in L^2_{\ov B}(\R^n,\Wedge^{\ell})$; and
\item[(iii)] if $B\in\mathscr{B}$ then there exists an $n$-form $\vartheta^B$ of $\Cinf_0(\R^n,\Wedge^n)$ supported in $\frac{1}{2}\,B$ such that $dT^B_nu=u-(\int u)\vartheta^B$ for every $u\in L^2_{\ov B}(\R^n,\Wedge^n)$.
\end{itemize}

We return now to the proof of the lemma. Suppose that $a$ is an $H^p_d(\R^n,\Wedge^{\ell})$-atom, where $1\leq\ell\leq n-1$. Then there is a ball $B$ and a form $b'$ in $L^2(\R^n,\Wedge^{\ell-1})$ such that $\supp b'\subset \ov B$, $a=db'$ and $\norm{a}_2\leq|B|^{1/p-1/2}$. Set $b=T^B_{\ell}a$, noting that $\supp b\subset \ov B$. Moreover,
\[a=dT^B_{\ell}a+T^B_{\ell+1}da=db+T^B_{\ell+1}d^2b'=db\]
and
\[\norm{b}_{L^2(\R^n,\Wedge)}=\norm{T^B_{\ell}a}_{L^2(\R^n,\Wedge)} \leq c_nr(B)\norm{a}_{L^2(\R^n,\Wedge)} \leq c_nr(B)|B|^{1/2-1/p}.\]
This proves the lemma when $1\leq\ell\leq n-1$. The case when $\ell=n$ may be proved similarly.
\end{proof1}

Henceforth we suppose that $\Omega$ is a special Lipschitz domain with Lipschitz constant $A$. Recall that $H^p(\R^n)\subset\mathcal{S}'(\R^n)$, $\Fdot^0_{p,2}(\R^n)\subset\mathcal{Z}'(\R^n)$ and that the natural projection $\J$ from $\mathcal{S}'(\R^n)$ to $\mathcal{Z}'(\R^n)$ induces an isomorphism from $H^p(\R^n)$ to $\Fdot^0_{p,2}(\R^n)$. In this way, the theory already developed in this paper can be applied, because, when $0<p\leq1$, 
\begin{align*}
\J H^p(\R^n,\Wedge)&= \dot F^0_{p,2}(\R^n,\Wedge);\\ 
\J H^p_d(\R^n,\Wedge)&= d\dot F^1_{p,2}(\R^n,\Wedge); \quad \text{and}\\  
\J H^p_{z,d}(\ov\Omega,\Wedge)&= d\dot F^1_{p,2, \ov\Omega}(\R^n,\Wedge)
\end{align*}
with equivalence of norms. The second and third identities follow from Theorem \ref{th:T on special lipschitz} and Corollary \ref{cor:T on special lipschitz}. In particular, it is a consequence of Corollary~\ref{cor:T on special lipschitz} that $dT$, correctly interpreted, is the identity on $H^p_{z,d}(\ov\Omega,\Wedge)$.  Using the notation of Section \ref{sec:CKS}. including the condition \eqref{meanzero}, we obtain a Calder\'on-type reproducing formula on this space, namely
\begin{equation}\label{calderon}
\J u=dT\J u=\int_0^\infty d(\phi_t\ast\Psi_t\ast \contr \J u)\,dt
=\int_0^\infty(\partial\phi)_t\ast\wedg
\,(\Psi_t\ast \contr \J u)\,\frac{dt}{t}\,.
\end{equation} 
We note that this formula actually holds on the whole space
$H^p_{d}(\R^n,\Wedge)$.  By Remark \ref{convergence} and the fact that ${\mathcal Z}(\R^n)$ is dense in $\Fdot^0_{p,2}(\R^n)$,
 the implicit limits in these improper integrals exist in  $H^p(\R^n,\Wedge)$ whenever $u\in H^p(\R^n,\Wedge)$.

The next two theorems, which characterise the spaces $H^p_d(\R^n,\Wedge)$ and $H^p_{z,d}(\ov\Omega,\Wedge)$ in terms of atoms, are the two main results of this section.

\begin{theorem}\label{th:atomic characterisation of  H^p_d}
Suppose that $1\leq\ell\leq n$ and $n/(n+1)<p\leq1$. There exist constants $c_p$ and $c_p'$ with the following properties.
\begin{enumerate}
\item[(i)] If $(a_k)_{k=0}^{\infty}$ is a sequence of $H^p_d(\R^n,\Wedge^{\ell})$-atoms and $(\lambda_k)_{k=0}^{\infty}$ belongs to $\ell^p(\C)$ then the series
\[\sum_{k=0}^{\infty}\lambda_ka_k\]
converges in $H^p(\R^n,\Wedge)$ to a form $u$ in $H^p_d(\R^n,\Wedge^{\ell})$, and
\begin{equation}\label{eq:Hp bounded above}
\norm{u}_{H^p(\R^n,\Wedge)}^p\leq c_p\sum_{k=0}^{\infty}|\lambda_k|^p.
\end{equation}
\item[(ii)] Conversely, if $u\in H^p_d(\R^n,\Wedge^{\ell})$ then there is a sequence $(a_k)_{k=0}^{\infty}$ of $H^p_d(\R^n,\Wedge^{\ell})$-atoms and a sequence $(\lambda_k)_{k=0}^{\infty}$ in $\ell^p(\C)$ such that
\[u=\sum_{k=0}^{\infty}\lambda_ka_k,\]
where the sum converges in $H^p(\R^n,\Wedge)$, and
\begin{equation}\label{eq:Hp bounded below}
\sum_{k=0}^{\infty}|\lambda_k|^p\leq c_p'\norm{u}_{H^p(\R^n,\Wedge)}^p.
\end{equation}
\end{enumerate}
\end{theorem}

\begin{theorem}\label{th:atomic characterisation of  H^p_zd}
Suppose that $1\leq\ell\leq n$, $n/(n+1)<p\leq1$ and $\Omega$ is a special Lipschitz domain in $\R^n$. Then there   exist constants $c_p$ and $c_p'$ with the following properties.
\begin{enumerate}
\item[(i)] If $(a_k)_{k=0}^{\infty}$ is a sequence of $H^p_{z,d}(\ov\Omega,\Wedge^{\ell})$-atoms and $(\lambda_k)_{k=0}^{\infty}$ belongs to $\ell^p(\C)$ then the series
\[\sum_{k=0}^{\infty}\lambda_ka_k\]
converges in $H^p(\R^n,\Wedge)$ to a form $u$ in $H^p_{z,d}(\ov\Omega,\Wedge^{\ell})$, and
\[\norm{u}_{H^p(\R^n,\Wedge)}^p\leq c_p\sum_{k=0}^{\infty}|\lambda_k|^p.\]
\item[(ii)] Conversely, if $u\in H^p_{z,d}(\ov\Omega,\Wedge^{\ell})$ then there is a sequence $(a_k)_{k=0}^{\infty}$ of $H^p_{z,d}(\ov\Omega,\Wedge^{\ell})$-atoms and a sequence $(\lambda_k)_{k=0}^{\infty}$ in $\ell^p(\C)$ such that
\[u=\sum_{k=0}^{\infty}\lambda_ka_k,\]
where the sum converges in $H^p(\R^n,\Wedge)$, and
\[\sum_{k=0}^{\infty}|\lambda_k|^p\leq c_p'\norm{u}_{H^p(\R^n,\Wedge)}^p.\]
\end{enumerate}
\end{theorem}

The results of the preceding theorems are generalisations to exact forms of the classical atomic decompositions of \cite{CMS85} for $\R^n$ and \cite{CKS93} for special Lipschitz domains $\Omega\subset\R^n$. The generalisation to exact forms first appeared in \cite{LM05} and \cite{LM04}  for the special case when $p=1$. Apart from expanding the range of $p$, our contribution is a new proof using the reproducing formula (\ref{calderon}), which is especially suited for application to special Lipschitz domains due to the support properties of $\phi$. Consequently, our proof of Theorem~\ref{th:atomic characterisation of  H^p_zd} is shorter and more direct than the one given in \cite{LM04}, since we avoid using reflection maps and obtain more efficiently the desired support properties for $H^p_{z,d}(\ov\Omega,\Wedge^{\ell})$-atoms. As a by-product of our proof, we also obtain a special atomic decomposition for tent space functions supported in tents over $\Omega$ (see Theorem~\ref{th:tent space decomposition}).

Before we can prove these characterisations, it is necessary to present a sequence of definitions and lemmata related to `tents' over open sets, `tent spaces' and tent space atoms. To help the reader contextualise what follows, we first offer a brief outline of the proof of Part (ii) of each of the above theorems. Suppose that $u\in H^p_d(\R^n,\Wedge^{\ell})$. Following the method developed in \cite{CMS85}, we define an operator $Q$ by
\[(Qu)(x,t)={\Psi}_t\ast\contr u(x)\qquad  \forall\,(x,t)\in\R^n\times\R^+=\R^n\times(0,\infty)\, ,\] 
and show that $Qu$ belongs to the tent space $\tent^p(\R^n\times\R^+,\Wedge^{\ell-1})$. Using the classical atomic decomposition for tent space functions, one may write $Qu=\sum_k\lambda_kA_k$, where each $A_k$ is a $\tent^p(\R^n\times\R^+,\Wedge^{\ell-1})$-atom and the sequence $(\lambda_k)$ belongs to $\ell^p$. One then constructs a map $\pi$  by 
\[\pi U=\int_0^{\infty}(\partial\phi)_t\ast\wedg\, U(\cdot,t)\,\frac{dt}{t}\]
so that, by the reproducing formula (\ref{calderon}),
\[u = dTu= \pi Qu = \sum_k\lambda_k\pi A_k\]
where each $\pi A_k$ is an $H^p_d(\R^n,\Wedge^{\ell})$-atom,
thus  obtaining the atomic decomposition for $u$.

The atomic decomposition for elements of $H^p_{z,d}(\ov\Omega,\Wedge^{\ell})$ will be proved along the same lines with the following variations. If $u\in H^p_{z,d}(\ov\Omega,\Wedge^{\ell})$ then $Qu$ is in fact supported in a `tent' over $\Omega$. So the tent space decomposition for $Qu$ (presented in Theorem \ref{th:tent space decomposition}) gives tent space atoms $A_k$ with good support properties with respect to the domain $\Omega$. It follows that $u=\pi Qu=\sum_k\lambda_k\pi A_k$, where each $\pi A_k$ can be written as a finite sum of $H^p_{z,d}(\ov\Omega,\Wedge^{\ell})$-atoms.

We turn now to the relevant definitions. If $\app>0$ and $x\in\R^n$, let $\Gamma'_{\app}(x)$ denote the cone in $\R^n\times\R^+$ with aperture $\app$ and vertex at $x$, namely
\[\Gamma'_{\app}(x)=\{(y,t)\in\R^n\times\R^+:|y-x|<\app t\}.\]
If $O$ is an open subset of $\R^n$, then the tent $T_{\app}(O)$ over $O$ with aperture $\app$ is defined by
\[T_{\app}(O)=\{(y,t)\in\R^n\times\R^+:d(y,O^c)\geq\app t\}.\]
We follow the convention of writing $\Gamma'(x)$ for $\Gamma_1'(x)$ and  $T(O)$ for $T_1(O)$. Given any measurable function $U$ on $\R^n\times\R^+$, we define the Lusin area integral $SU$ of $U$ by the formula
\[(SU)(x)=\left(\iint_{\Gamma'(x)}|U(y,t)|^2\frac{dy\,dt}{t^{n+1}}\right)^{1/2}.\]

\begin{definition}
Suppose that $p>0$. The tent space $\tent^p(\R^n\times\R^+)$ is defined to be the set of all measurable functions $U$ on $\R^n\times\R^+$ such that $\norm{U}_{\tent^p(\R^n\times\R^+)}$ is finite, where
\[\norm{U}_{\tent^p(\R^n\times\R^+)}=\norm{SU}_{L^p(\R^n)}.\]
\end{definition}

The tent spaces were first introduced in the article \cite{CMS85} of Coifman, Meyer and Stein, and owe their name to the fact that, when $0<p\leq1$, their functions can be decomposed as a sum of atoms supported in tents over balls.

\begin{definition}\label{def:T^p atom} Suppose that $p>0$. A measurable function $A$ on $\R^n\times\R^+$ is said to be a $\tent^p(\R^n\times\R^+)$-atom if there exists a ball $B$ in $\R^n$ such that $\supp A\subset T(B)$ and
\[\iint_{\R^n\times\R^+}|A(y,t)|^2\, dy\,\frac{dt}{t}\leq |B|^{1-2/p}.\]
\end{definition}

If $0<p\leq1$ then it is relatively straightforward to show that every $\tent^p(\R^n\times\R^+)$-atom $A$ belongs to $\tent^p(\R^n\times\R^+)$ and that $\norm{A}_{\tent^p(\R^n\times\R^+)}\leq1$. Consequently, if $(\lambda_k)_{k\in\N}\in\ell^p(\C)$ and $(A_k)_{k\in\N}$ is a sequence of $\tent^p(\R^n\times\R^+)$-atoms then $\sum_{k\in\N}\lambda_kA_k$ belongs to $\tent^p(\R^n\times\R^+)$. That the following converse is true is a deeper result due to Coifman, Meyer and Stein \cite{CMS85}.

\begin{theorem}\label{th:atomic decomposition of tent spaces}
Suppose that $0<p\leq1$. There exists a constant $C$ (depending only on $n$ and $p$) with the following property: for all $U$ in $\tent^p(\R^n\times\R^+)$, there exists a sequence $(\lambda_k)_{k\in\N}$ in $\ell^p(\C)$ and a sequence $(A_k)_{k\in\N}$ of $\tent^p(\R^n\times\R^+)$-atoms such that
\[U=\sum_{k\in\N}\lambda_kA_k\]
and
\[\sum_{k\in\N}|\lambda_k|^p\leq C\norm{U}^p_{\tent^p(\R^n\times\R^+)}.\]
\end{theorem}

We introduce the following variant of the above theorem, where the tent space atoms are supported in Carleson boxes with good support properties with respect to an underlying domain $\Omega$.

\begin{theorem}\label{th:tent space decomposition}
Suppose that $0<p\leq1$ and $\beta>0$. There exist positive constants $C'$ (depending only on $n$, $p$ and $\beta$) and $c_{\beta}$ (depending only on $\beta$), where $0<c_{\beta}<1$, satisfying the following property. If $U\in\tent^p(\R^n\times\R^+)$ and $\supp U\subset T_{\beta}(\Omega)$ for some proper open subset $\Omega$ of $\R^n$, then there exists a sequence $(\lambda_k)_{k\in\N}$ in $\ell^p(\C)$ and a sequence $(A_k)_{k\in\N}$ of $\tent^p(\R^n\times\R^+)$-atoms, supported in corresponding  tents $(T(B_k))_{k\in\N}$,  such that
\begin{enumerate}
\item[(i)] $\ds U=\sum_{k\in\N}\lambda_kA_k$,
\item[(ii)] $\ds \sum_{k\in\N}|\lambda_k|^p\leq C'\norm{U}^p_{\tent^p(\R^n\times\R^+)}$, 
\item[(iii)] $\supp A_k\subset c_{\beta}B_k\times(0,6\beta^{-1}c_{\beta}\,r(B_k))$ and  $5c_{\beta}B_k\subset\Omega$ whenever $k\in\N$.
\end{enumerate}
\end{theorem}

\begin{proof1}
The proof is an adaptation of the proof of \cite[Theorem 1.1]{eR07}, which in turn is based on the original ideas presented in \cite{CMS85}. Fix any number $\nu$ in the interval $(0,1)$.  For any $k$ in $\Z$, let $O_k$ denote the open subset of $\R^n$ given by
\[O_k=\{x\in\R^n:SU(x)>2^k\}.\]
It can be shown that
\begin{equation}\label{eq:covering of supp U}
\supp U\subset \bigcup_{k\in\Z}T_{\nu}(O_k^*),
\end{equation}
where each open set $O_k^*$ is constructed using a corresponding set of global $\gamma$-density (see \cite[pp.~128--130]{eR07} for details). For each integer $k$, the Whitney lemma (see, e.g., \cite[Lemma 2.2]{eR07}) applied to the open set $O_k^*\cap\Omega$ gives a denumerable index set $I^k$, a sequence of balls $(B_j^k)_{j\in I^k}$ having radii  $(r_j^k)_{j\in I^k}$ and centres $(x_j^k)_{j\in I^k}$, and a sequence $(\varphi_j^k)_{j\in I^k}$ of nonnegative functions on $\R^n$ with the following properties:
\[
O_k^*\cap\Omega=\bigcup_{j\in I^k}B_j^k,\ \ 
d(x_j^k,(O_k^*\cap\Omega)^c)=10r_j^k,\ \ 
\supp\varphi_j^k\subset 2B_j^k,\ \ 
\sum_{j\in I^k}\varphi_j^k=\ind_{O^*_k\cap\Omega}\,,
\]
and
\begin{equation}\label{disjoint}
\tfrac{1}{4}B_i^k\cap\tfrac{1}{4}B_j^k=\emptyset \qquad\text{if }i\neq j.
\end{equation}
It can be shown that $O^*_{k+1}\subset O^*_k$ for all $k$ (see \cite[pp. 128, 130]{eR07}). Therefore, for each $(x,t)$ in $\R^n\times\R^+$,
\[
\left(\ind_{T_{\nu}(O^*_k)}-\ind_{T_{\nu}(O^*_{k+1})}\right)(x,t)\,\ind_{\Omega}(x)
=\sum_{j\in I^k}\varphi^k_j(x)\left(\ind_{T_{\nu}(O^*_k)}-\ind_{T_{\nu}(O^*_{k+1})}\right)(x,t)\,\ind_{\Omega}(x)
\]
and hence
\[U(x,t)=\sum_{k\in\Z}\sum_{j\in I^k}U(x,t)\varphi^k_j(x)\left(\ind_{T_{\nu}(O^*_k)}-\ind_{T_{\nu}(O^*_{k+1})}\right)(x,t)\]
by (\ref{eq:covering of supp U}). Define, for all integers $k$ and all $j$ in $I^k$,
\begin{align*}
\mu_j^k&=\int_0^{\infty}\int_{\R^n}|U(y,t)|^2\varphi_j^k(y)^2\left(\ind_{T_{\nu}(O^*_k)}-\ind_{T_{\nu}(O^*_{k+1})}\right)(y,t)\,dy\,\frac{dt}{t}\,,\\
A^k_j(y,t)&=U(y,t)\varphi_j^k(y)\left(\ind_{T_{\nu}(O^*_k)}-\ind_{T_{\nu}(O^*_{k+1})}\right)(y,t)|B_j^k|^{1/2-1/p}(\mu^k_j)^{-1/2} 
\intertext{(unless $\mu_j^k=0$, in which case we define $A_j^k=0$) and}		
\lambda_j^k&=|B_j^k|^{1/p-1/2}(\mu^k_j)^{1/2}\,.
\end{align*}
Then
\[U=\sum_{k\in\Z}\sum_{j\in I^k}\lambda_j^kA_j^k.\]

We claim that, up to a multiplicative constant, each $A_j^k$ is a $\tent^p(\R^n\times\R^+)$-atom with the desired properties. 

First we show that $\supp A_j^k\subset T(\tilde{B}_j^k)$, where $\tilde{B}_j^k=CB_j^k$ and $C=2+12/\max\{\beta,\nu\}$. Suppose that $(y,t)\in T_{\nu}(O^*_k)\cap T_{\beta}(\Omega)$ and $y\in\supp\varphi_j^k$; that is,
\[d(y,(O^*_k)^c)\geq\nu t, \quad d(y,\Omega^c)\geq\beta t \quad\mbox{and}\quad |y-x_j^k|<2r^k_j.\] 
We aim to show that $d(y,(CB_j^k)^c)\geq t$ for then $\supp A_j^k\subset T(\tilde{B}_j^k)$. Suppose that $z\in(CB_j^k)^c$. Then
\begin{equation}\label{eq:|y-z| estimate}
|y-z|\geq |z-x_j^k|-|y-x_j^k|\geq (C-2)r_j^k=\frac{12r_j^k}{\max\{\beta,\nu\}}.
\end{equation}
Also, $d(x_j^k,(O^*_k\cap\Omega)^c)=10r_j^k$. Suppose that $\varepsilon>0$. There exists $u$ in $(O^*_k\cap\Omega)^c$ such that $|x_j^k-u|<10r_j^k+\varepsilon$. So
\[\max\{\beta,\nu\}t \leq |y-u| \leq |y-x_j^k|+|x_j^k-u|<12r_j^k+\varepsilon.\]
Since this is true for every positive $\varepsilon$, it follows that $\max\{\beta,\nu\}t \leq 12r_j^k$. Combining this with (\ref{eq:|y-z| estimate}) gives $|y-z|\geq t$, and hence $d(y,(CB_j^k)^c)\geq t$ as required.

Second, the definition of $A_j^k$ implies that
\[\iint|A_j^k(y,t)|^2\,dy\,\frac{dt}{t}=|B_j^k|^{1-2/p}=C^{n(2/p-1)}|\tilde{B}_j^k|^{1-2/p},\]
and so up to the multiplicative constant $C^{n(2/p-1)}$, each $A_j^k$ is a $\tent^p(\R^n\times\R^+)$-atom.

Third, we prove  that the $A^k_j$ satisfy support properties as in  Part (iii) of the theorem. Now each $A_j^k$ is supported in $T_{\beta}(\Omega)\cap (2B_j^k\times\R^+)$, where $5(2B_j^k)\subset\Omega$. So if $(y,t)\in T_{\beta}(\Omega)\cap  (2B_j^k\times\R^+) $ then
\[\beta t\leq d(y,\Omega^c)\leq d(y,x_j^k)+d(x_j^k,\Omega^c)< 2r_j^k+10r_j^k\]
and hence $0<t<12\beta^{-1}r_j^k$. This shows that $\supp A_j^k\subset 2B_j^k\times(0,12\beta^{-1}r_j^k)$. Defining the constant $c_{\beta}$ by $c_{\beta}=2/C$, it is now easy to see that
\[\supp A_j^k\subset T(\tilde{B}_j^k), \quad  \supp A^k_j\subset c_{\beta}\tilde{B}_j^k\times(0,6\beta^{-1}c_{\beta}\,r(\tilde{B}_j^k))   \quad\text{and}\quad 5c_{\beta}\tilde{B_j^k}\subset\Omega\,.\]
 It  remains to show that there exists a constant $C'$, independent of $\Omega$ and $U$, such that
\[\sum_{k\in\Z}\sum_{j\in I^k}|\lambda_j^k|^p\leq C'\norm{U}_{\tent^p(\R^n\times\R^+)}\,.\]
The proof, which uses (\ref{disjoint}), proceeds exactly as in \cite[pp.~132--133]{eR07} and will not be reproduced here. 
 Finally, relabel the balls $(\tilde B^k_j)$ as $(B_k)$, $(C^{n(1/p-1/2)}\lambda^k_j)$ as $(\lambda_k)$ and the functions $(C^{n(1/2-1/p)}A^k_j)$ as $(A_k)$.
This completes the proof of the theorem.
\end{proof1}

\begin{remark}
A comparison of this proof with the proof of  \cite[Theorem 1.1]{eR07} shows that Theorem~\ref{th:tent space decomposition} also holds when the underlying space $\R^n$ is replaced by any space $X$ of homogeneous type that satisfies the assumptions of \cite[\S 1.1.3]{eS93}.
\end{remark}

 \begin{remark} Theorem \ref{1.1} follows from the special case of Theorem \ref{th:tent space decomposition} when $p=1$ and $\beta=1$, with $(c_\beta B_k), (c_\beta^{n/2}\lambda_k)$ and $(c_\beta^{-n/2}A_k)$ renamed as $(B_k), (\lambda_k)$ and $(A_k)$. 
\end{remark}

Recall that $A$ and $\ap$ are fixed positive numbers such that $\ap A<1$. In order to apply the previous theorem, we need the following.

\begin{lemma}\label{lem:supp Qu}
Suppose that  $a>0$, $\Omega$ is a special Lipschitz domain with  Lipschitz constant $A$ and that $\psi$ is a $\Cinf(\R^n)$ function supported in $\{y\in\overline{\Gamma_{\ap}(0)}:\,y_n\geq a\}$.  Suppose also that $u\in\mathcal{S}'(\R^n)$ and $\supp u\subset\ov\Omega$. Define $Qu$ by 
\[Qu(x,t)=(\psi_t\ast u)(x)\qquad\forall(x,t)\in\R^n\times\R^+.\]
Then $Qu$ is supported in $T_{\beta}(\Omega)$ where
\begin{equation}\label{eq:beta}
\beta=\frac{a(1-\ap A)}{\sqrt{1+A^2}}\ .
\end{equation} 
\end{lemma}

\begin{proof1}
What needs to be shown is that, under the stated hypotheses, $\text{dist}(\supp(\psi_t*u),\Omega^c)\geq\beta t$. By (\ref{sppt}), $\supp(\psi_t\ast u)\subset \overline\Omega+\{y\in\overline{\Gamma_{\ap}(0)}:y_n\geq at\}$, so we need to show that if $x\in\overline\Omega$, $y\in\overline{\Gamma_\ap(0)}$, $y_n\geq at$ and $z\in\Omega^c$, then $|x+y-z|\geq \beta t$. Let $w=z-x$, and note that by the assumption on $\Omega$, $w_n\leq A |w'|$.

So the result is proved once we show that $|y-w|\geq \beta t$ whenever $-\infty<w_n\leq A|w'|$, $|y'|\leq \ap y_n$ and $y_n\geq at$. We split into two cases.
\newline Case (i): $w_n\leq \ap A y_n$. Then 
\[|y-w|\geq y_n-w_n \geq (1-\ap A)y_n\geq (1-\ap A)at > \beta t\]
by (\ref{eq:beta}).
\newline Case (ii): $w_n> \ap A y_n$. Then $|w'|\geq \frac1A w_n>\ap y_n\geq |y'|$. So
\begin{align*}
|y-w|^2 & = (y_n-w_n)^2 +|y'-w'|^2\\
&> (y_n-w_n)^2 + (|w'|-|y'|)^2\\
&> (y_n-w_n)^2 + (\tfrac1A w_n-\ap y_n)^2\\
&=\tfrac{1+A^2}{A^2}w_n^2 -\tfrac{2(\ap+A)}A y_nw_n
+(1+\ap^2)y_n^2\\
&\geq \left(\tfrac{\beta}{a}\right)^2 y_n^2 \\
&\geq \beta^2t^2\,,
\end{align*}
where we have minimised over $w_n$ in the usual way for quadratic expressions.
\end{proof1}

We are now in a position to prove Theorems~\ref{th:atomic characterisation of  H^p_d} and \ref{th:atomic characterisation of  H^p_zd}.

\begin{proof1}
Suppose throughout that $n/(n+1)<p\leq1$.

First we prove Part (i) of Theorem~\ref{th:atomic characterisation of  H^p_d}. 
Suppose that $(a_k)_{k=0}^{\infty}$ is a sequence of $H^p_d(\R^n,\Wedge)$-atoms, $(\lambda_k)_{k=0}^{\infty}$ belongs to $\ell^p(\C)$ and $a_k=db_k$. Since each component of $a_k$ is a classical $H^p(\R^n)$-atom, the classical theory implies that there exist a constant $c_p$ and $u$ in $H^p(\R^n,\Wedge)$ such that (\ref{eq:Hp bounded above}) holds and
\begin{equation}\label{eq:u sum}
u=\sum_{k=0}^{\infty}\lambda_ka_k=\sum_{k=0}^{\infty}\lambda_kdb_k,
\end{equation}
where the sum converges in $H^p(\R^n,\Wedge)$.

Recall that $\J$ denotes the natural projection from $\mathcal{S}'(\R^n)$ to $\mathcal{Z}'(\R^n)$. By (\ref{eq:u sum}),
\[\sum_{k=0}^M \lambda_k\J a_k=d\Big(\sum_{k=0}^M \lambda_k\J b_k\Big)\to \J u \qquad\text{in $\Fdot^0_{p,2}(\R^n,\Wedge)$ as $M\to\infty$},\]
and it follows from the continuity of $d$ and Theorem~\ref{th:T on special lipschitz} that $\J u\in d\Fdot^1_{p,2}(\R^n,\Wedge)=\J H^p_d(\R^n,\Wedge)$. Hence $u\in H^p_d(\R^n,\Wedge)$. This completes the proof of Theorem~\ref{th:atomic characterisation of  H^p_d} (i).

Part (i) of Theorem~\ref{th:atomic characterisation of  H^p_zd} is proved along the same lines, noting that $d\Fdot^1_{p,2,\ov\Omega}(\R^n,\Wedge)=\J H^p_d(\ov\Omega,\Wedge)$.

We now prove the converse statements in each of Theorems~\ref{th:atomic characterisation of  H^p_d} and \ref{th:atomic characterisation of  H^p_zd}. Let $\phi$ denote the $\Cinf_0(\R^n)$ function of Proposition \ref{phi}, chosen such that
\begin{equation}\label{eq:momentcondition1}\int_{\R^n}{\Psi}(x)\,dx=0
\end{equation}
where $\Psi(x)=2\phi(x)x$.  We remark also that $\int_{\R^n}\phi(x)\,dx=1$, and
\begin{equation}\label{eq:momentcondition2}
 \qquad  \int_{\R^n}(\partial\phi)(x)\,dx=0 
\end{equation}
where $\partial \phi=\Sigma_{j=1}^n\partial_j\phi\, e_j$,

Given $u$ in $\mathcal{S}'(\R^n,\Wedge)$, define $Qu$ by
\[(Qu)(x,t)={\Psi}_t\ast\contr\, u(x)\]
whenever $(x,t)\in\R^n\times\R^+$. By the moment condition 
(\ref{eq:momentcondition1}), it is well-known (see, for example, \cite[p.~308]{CKS93}) that $Q$ is bounded from $H^p(\R^n,\Wedge)$ to $\tent^p(\R^n\times\R^+,\Wedge)$. Given $U$ in $\tent^2(\R^n\times\R^+,\Wedge)$ with compact support in $\R^n\times\R^+$, define $\pi U$ by the formula
\[\pi U=\int_0^{\infty}(\partial\phi)_t\ast\wedg\, U(\cdot,t)\,\frac{dt}{t}.\]
Again, by the moment condition (\ref{eq:momentcondition2}), it is well known (see \cite[Theorem 6]{CMS85}) that the operator $\pi$ extends to a bounded linear operator from $\tent^2(\R^n\times\R^+,\Wedge)$ to $L^2(\R^n,\Wedge)$ and from $\tent^p(\R^n\times\R^+,\Wedge)$ to $H^p(\R^n,\Wedge)$.


We focus now on the proof of Theorem~\ref{th:atomic characterisation of  H^p_d}~(ii). Suppose that $u\in H^p_d(\R^n,\Wedge^{\ell})$, where $1\leq\ell\leq n$ and $n\geq1$. Then $Qu$ belongs to the tent space $\tent^p(\R^n\times\R^+,\Wedge^{\ell-1})$ and by Theorem~\ref{th:atomic decomposition of tent spaces}, $Qu$ has the atomic decomposition
\[Qu=\sum_{k\in\Z}\lambda_kA_k,\]
where each $A_k$ is a $\tent^p(\R^n\times\R^+,\Wedge^{\ell-1})$-atom supported in a tent $T(B_k)$, the sum converges in $\tent^p(\R^n\times\R^+,\Wedge)$ and
\[\sum_k|\lambda_k|^p<c_p''\norm{Qu}_{\tent^p(\R^n\times\R^+,\Wedge)}^p \leq c'_p\|u\|_{H^p(\R^n,\Lambda)}\,.\]
Define $a_k$ by $a_k=\pi A_k$, so that
\begin{align*}a_k(x)
&=\int_0^{\infty}\int_{\R^n}(\partial\phi)_t(x-y)\wedg A_k(y,t)\,dy\,\frac{dt}{t}\\
&=d\int_0^{\infty}\int_{\R^n}\phi_t(x-y)A_k(y,t)\,dy\,dt\\
&=db_k(x),
\intertext{where}
b_k&= \int_0^{\infty}\phi_t\ast A_k(\,\cdot\,,t)\,dt.
\end{align*}
We claim that each $a_k$ is an $H^p_d(\R^n,\Wedge^{\ell})$-atom, up to a multiplicative constant independent of $k$. First,
\[\norm{a_k}_{L^2(\R^n,\Wedge)}= \norm{\pi A_k}_{L^2(\R^n,\Wedge)}\leq C\norm{A_k}_{\tent^2(\R^n\times\R^+,\Wedge)}\leq C|B_k|^{1/2-1/p},\]
where the constant $C$ is independent of $k$. Second, we show that $b_k\in L^2(\R^n,\Wedge)$. Suppose that $B_k=B_{r_k}(z_k)$.  Successive applications of the triangle, Cauchy--Schwarz and Young's inequalities yield
\begin{align*}
\norm{b_k}_{L^2(\R^n,\Wedge)}^2
& \leq r_k \int_0^{r_k}\norm{\phi_t\ast A_k(\,\cdot\,,t)}_{L^2(\R^n,\Wedge)}^2  dt \\
& \leq r_k\int_0^{r_k}\norm{\phi_t}_{L^1(\R^n)}^2\norm{A_k(\,\cdot\,,t)}_{L^2(\R^n,\Wedge)}^2\,dt\\
&= C_1^2 r_k \int_0^{r_k} t\int_{\R^n}|A_k(y,t)|^2\,dy\,\frac{dt}{t}\\
&\leq C_1^2 r_k^2\,|B_k|^{1-2/p},
\end{align*} 
where $C_1=\|\phi\|_1$, and the final estimate follows from the fact that $A_k$ is a $\tent^p(\R^n\times\R^+,\Wedge^{\ell})$-atom. Third, we note that $\supp a_k\subset B_k$. Indeed, by (\ref{sppt}) and (\ref{set}),
\begin{align*} 
\supp a_k=
\supp(\pi A_k) &\subset\bigcup_{0\leq t\leq r_k}\left\{\supp \phi_t+\left(1-\frac t {r_k}\right)B_k\right\}\\
&\subset\bigcup_{0\leq t\leq r_k}\left\{t B_1(0)+(r_k- t)B_1(z_k)\right\}\\
&=B_{r_k}(z_k)\\
&=B_k\,.
\end{align*}
Hence, up to a multiplicative constant, each $a_k$ is an $H^p_d(\R^n,\Wedge^{\ell})$-atom as claimed.

It remains to be shown that $u=\sum_k\lambda_k a_k$,  where the sum converges in the topology of $H^p(\R^n,\Wedge)$. Since $\sum_k\lambda_k A_k$ converges in $\tent^p(\R^n\times\R^+,\Wedge)$ and $\pi$ is bounded from $\tent^p(\R^n\times\R^+,\Wedge)$ to $H^p(\R^n,\Wedge)$, it follows that $\pi Qu=\sum_k\lambda_k a_k$, where the sum converges in the topology of $H^p(\R^n,\Wedge)$. But note by the definitions of $Q$ and $\pi$  that $\J\pi Qu=dT\J u$, where $T$ is the operator given by \eqref{reproducing}. Recall from (\ref{calderon}) that $dT$ is the identity on $d\Fdot^1_{p,2}(\R^n,\Wedge)$. Since $\J u\in d\Fdot^1_{p,2}(\R^n,\Wedge)$, we have that
\[\J u=dT\J u=\J\pi Qu=\J\sum_k\lambda_k a_k\]
and hence that $u=\sum_k\lambda_ka_k$ as required. 
This, together with the bound already proved on $\sum|\lambda_k|^p$, completes the proof of Theorem~\ref{th:atomic characterisation of  H^p_d} (ii).

We turn now to prove Theorem~\ref{th:atomic characterisation of  H^p_zd} (ii). Suppose that $u\in H^p_{z,d}(\ov\Omega,\Wedge^{\ell})$. Using the same argument as above, $u=\pi Qu=\sum_k\lambda_k\pi A_k$ where the sum converges in $H^p$, $\pi A_k=db_k$, $\supp A_k\subset T(B_k)$ and $\{\lambda_k\}\in\ell^p$. In this case, each $\pi A_k$ is not (even up to a multiplicative constant) necessarily an $H^1_{d,\ov\Omega}(\R^n,\Wedge^{\ell})$-atom because it may not satisfy the required support properties. We will instead show that each $\pi A_k$ can be written as a finite sum of $H^p_{z,d}(\ov\Omega,\Wedge^{\ell})$-atoms.

Since $\supp u\subset\ov\Omega$, we conclude by Lemma~\ref{lem:supp Qu} that $\supp Qu\subset \ov{T_{\beta}(\Omega)}$, where $\beta$ is given by (\ref{eq:beta})  with $a=1/2$. Theorem~\ref{th:tent space decomposition} gives the additional information that $A_k$ can be chosen so that $\supp (A_k)\subset  cB_k\times(0,6\beta^{-1}cr_k)$, 
and $5cB_k\subset\Omega$, where $c$ is independent of $k$. Hence $\text{dist}(\ov{cB_k},\Omega^c)\geq4cr_k$. By (\ref{sppt}) and (\ref{set}),
\begin{align*}
\supp a_k &=\supp(\pi A_k)\\
& \subset \bigcup_{0<t\leq 6\beta^{-1}cr_k}\left(\supp(\partial\phi)_t + \ov{cB_k}\,\right)\\
&\subset \ov{cB_k}+\{y\in\ov{\Gamma_\ap(0)}:\,y_n\leq 6\beta^{-1}cr_k\} =:G_k\,.
\end{align*}
Note that $G_k$ is a compact subset of $\Omega$ and that $\text{dist}(G_k,\Omega^c)\geq 4cr_k$. This is because, if $x\in\ov{cB_k}$, $y\in \Gamma_{\ap}(0)$ and $w\in \Omega^c$, then $w-y\in \Omega^c$, so $|(x+y)-w| = |x - (w-y)|\geq 4cr_k$.
So we may cover $G_k$ with finitely many balls $\{\tfrac{1}{2}B^j_k\}_{j=1}^M$ of radius $cr_k/2$ and centres $z_k^j$, where $z_k^j\in G_k$ and where (by scale and translation invariance) the integer $M$ is independent of $k$. Let $\{\eta^j_k\}_{j=1}^M$ denote a smooth subordinate partition of unity with the properties that $0\leq\eta_k^j\leq1$, $\supp\eta_k^j\subset B_k^j$,
\[\sum_{j=1}^M\eta_k^j(x)=1 \qquad\forall x\in G_k\]
and $\Vert\nabla\eta_k^j\Vert_{L^{\infty}(\R^n)}\leq c'r_k^{-1}$ for some constant $c'$ independent of $j$ and $k$. For each $k$ and $j$, define the function $a_k^j$ and scalar $\mu_k^j$ by
\[
a_k^j = (\mu_k^j)^{-1} d(\eta_k^jb_k) \qquad\mbox{and}\qquad
\mu_k^j = |B_k^j|^{1/p-1/2}\Vert d(\eta_k^jb_k)\Vert_{L^2(\R^n,\Wedge)}
 \]
Now each $a_k^j$ is an $H^p_{z,d}(\ov\Omega,\Wedge^{\ell})$-atom because $\eta_k^jb_k\in L^2(\R^n,\Wedge^{\ell-1})$, $\supp a_k^j\subset B_k^j$ where $4B_k^j\subset\Omega$, and $\Vert a_k^j \Vert_{L^2(\R^n,\Wedge)}\leq |B_k^j|^{1/2-1/p}$.

Note now that
\[u=\sum_{k=0}^{\infty}\lambda_k\pi A_k=\sum_{k=0}^{\infty}\sum_{j=1}^M\lambda_k\mu_k^ja_k^j,\]
where the sum converges in $H^p$ and $\sum_k|\lambda_k|^p<\infty$. To complete the proof of the theorem, it suffices to show that $\sup_{k,j}\mu_k^j\leq C'$ for some constant $C'$. This bound follows readily from the estimate
\begin{align*}
\Vert d(\eta_k^jb_k) \Vert_{L^2(\R^n,\Wedge)}
&=\norm{(d\eta_k^j)\wedg b_k+\eta_k^j\, db_k}_{L^2(\R^n,\Wedge)}\\
&\leq\Vert \nabla\eta_k^j\Vert_{L^{\infty}(\R^n)}\norm{b_k}_{L^2(\R^n,\Wedge)}+\Vert\eta_k^j\Vert_{L^{\infty}(\R^n)}\norm{\pi A_k}_{L^2(\R^n,\Wedge)}\\
&\leq c'r_k^{-1}r_k|B_k|^{1/2-1/p}+C|B_k|^{1/2-1/p}\\
&\leq (c'+C)|B_k|^{1/2-1/p}.
\end{align*}
Thus we have completed the proof of Theorem~\ref{th:atomic characterisation of  H^p_zd}.
\end{proof1}

\bigskip
\textbf{Addresses:}\\[1ex]
\textsc{Martin Costabel}\\
IRMAR, Universit\'e de Rennes 1,
Campus de Beaulieu\\35042 Rennes Cedex, France\\
E-mail: \texttt{martin.costabel@univ-rennes1.fr}
\\[1ex]
\textsc{Alan McIntosh}\\
Centre for Mathematics and its Applications, MSI\\
Australian National University\\
Canberra, ACT 0200, Australia\\
E-mail: \texttt{alan.mcintosh@anu.edu.au}
\\[1ex]
\textsc{Robert Taggart}\\
NSW Regional Forecasting Centre\\
Bureau of Meteorology\\
PO Box 413\\
Darlinghurst NSW 1300, Australia\\
E-mail: \texttt{r.taggart@bom.gov.au}

\end{document}